\documentclass[11pt]{article}

\usepackage{epsfig}

\usepackage{amsmath,makeidx,amssymb,amscd,amsfonts}
\usepackage{graphicx,epsfig,latexsym,bbm,index}
\usepackage[active]{srcltx}

\newcommand{\me} {1/2}

\newcommand{\tme} {3/2}

\newcommand{\eps}{\varepsilon}

\newcommand{\vphi}{\varphi}
\newcommand{\lbd}{\lambda}

\newtheorem{theorem}{Theorem}[section]
\newtheorem{lemma}[theorem]{Lemma}
\newtheorem{propo}[theorem]{Proposition}

\newtheorem{remark}[theorem]{Remark}

\def\proof{\par{\it Proof}. \ignorespaces}
\def\endproof{\vbox{\hrule height0.6pt\hbox{%
   \vrule height1.3ex width0.6pt\hskip0.8ex
   \vrule width0.6pt}\hrule height0.6pt }}

\def\N{{\mathord{\rm I\mkern-3.6mu N}}}

\def\R{{\mathord{\rm I\mkern-3.6mu R}}}

\begin{document}
\setcounter{footnote}{1}

\title{On level set type methods for elliptic Cauchy problems}

\author{
A.~Leit\~ao%
\thanks{Department of Mathematics, Federal University of St. Catarina,
P.O. Box 476, 88040-900 Florian\'opolis, Brazil (aleitao@mtm.ufsc.br)}
\ \ and \
M.~Marques Alves%
\thanks{IMPA, Est. D. Castorina 110, 22460-320 Rio de Janeiro, Brazil
(maicon@impa.br)} }

\date{\small\today}

\maketitle

\begin{small}
\noindent {\bf Abstract:}
Two methods of level set type are proposed for solving the Cauchy problem
for an elliptic equation.
Convergence and stability results for both methods are proven, characterizing
the iterative methods as regularization methods for this ill-posed problem.
Some numerical experiments are presented, showing the efficiency of our
approaches and verifying the convergence results.
\end{small}

\section{Introduction} \label{sec:1}

We start by introducing the inverse problem under consideration.
Let $\Omega \subset \R^d$, $d = 2, 3$, be an open bounded set with
piecewise Lipschitz boundary $\partial\Omega$. Further, we assume that
$\partial\Omega = \overline{\Gamma}_1 \cup \overline{\Gamma}_2$, where
$\Gamma_i$ are two open connected disjoint parts of $\partial\Omega$.
We denote by $P$ the elliptic operator defined in $\Omega$ by
\begin{equation} \label{def:P}
     P(u) \ := \ - \sum_{i,j=1}^{d}\,{D_i (a_{i,j} D_j u)} \, ,
\end{equation}
where the real functions $a_{i,j} \in L^\infty(\Omega)$ are such that the
matrix $A(x) := (a_{i,j})_{i,j=1}^d$ satisfies
$\xi^t A(x)\, \xi > \alpha ||\xi||^2$,
for all $\xi \in \R^d$ and for a.e. $x \in \Omega$, where $\alpha > 0$.

We denote by {\em elliptic Cauchy problem} the following boundary value
problem (BVP)\\[2ex]
$ (CP) \hskip3.5cm
  \left\{ \begin{array}{rl}
    P u   = f  \, ,&\!\!\! \mbox{in } \Omega \\
    u     = g_1\, ,&\!\!\! \mbox{at } \Gamma_1 \\
    u_\nu = g_2\, ,&\!\!\! \mbox{at } \Gamma_1
  \end{array} \right. . $ \\[2ex]
The functions $g_1, g_2: \Gamma_1 \to \R$ are given {\em Cauchy data},
and $f: \Omega \to \R$ is a known source term in the model.

As a solution of the elliptic Cauchy problem (CP) we consider every
$H^m(\Omega)$--distri\-bution, which solves the weak formulation of the
elliptic equation in $\Omega$ and also satisfies the Cauchy data at
$\Gamma_1$ in the sense of the trace operator ($m \in \N$ still has to
be chosen). Notice that, if we know either the Neumann or the Dirichlet
trace of $u$ at $\Gamma_2$, then $u$ can be computed as a solution of a
mixed BVP in a stable way. Therefore, it is enough to consider the task
of determining a trace (Dirichlet or Neumann) of $u$ at $\Gamma_2$.

It is well known that elliptic Cauchy problems are not well posed in the
sense of Hadamard.
A famous example given by Hadamard in the early 20's \cite{Had32, Le00}
shows that one cannot expect the solution of (CP) to depend continuously
on the data. For Lipschitz bounded domains $\Omega \subset \R^2$, the
severely ill-posedness of (CP) was recently investigated in \cite{Be07}
using a Steklov-Poincar\'e approach.
Existence of solutions for arbitrary Cauchy data $(g_1,g_2)$ cannot be
assured as a direct argumentation with the Schwartz reflection principle
shows \cite{GT77} (the Cauchy data $(g_1,g_2)$ are called {\em consistent}
if the corresponding problem (CP) has a solution).
It has been recently shown \cite{ABA06} that in the case $m = 1$, there
exists a dense subset $M$ of $H^{\me}(\Gamma_1) \times
[H^{\me}_{00}(\Gamma_1)]'$ such that (CP) has a $H^1(\Omega)$ solution
for Cauchy data $(g_1,g_2)\in M$.
What concerns uniqueness of solution for (CP), it is possible to extend the
Cauchy--Kowalewsky and Holmgren theorems to the $H^m$--context and prove
uniqueness of weak solutions (see, e.g., \cite{EL01} for the case $m = 1$).
For classical uniqueness results we refer the reader to \cite{Cal58}.
A weak uniqueness result for nonlinear Cauchy problems can be found in
\cite{KL03}.

A variety of numerical methods for solving (CP) can be found in the
literature.
An optimization approach based on least squares and Tikhonov regularization
was used in \cite{FM86}.
In \cite{KMF91} Mazya et al proposed an iterative algorithm based on the
successive solutions of well posed mixed BVPs.
A generalization of this method, based on fixed point theory, was derived in
\cite{Le00}.
A further generalization for nonlinear elliptic Cauchy problems
can be found in \cite{KL03}.
In \cite{Le98,HW01} the Backus-Gilbert method was used to solve (CP).
A Mann iterative regularization method was proposed for (CP) in \cite{EL01}.
In \cite{HL00} a method of conjugate gradient type was investigated.
Finite element approximations, based on an optimal control formulation of
(CP), are discussed in \cite{CN06}.
An application of the quasi-reversibility method for (CP) is considered in
\cite{Bou05,Bou06}.
In \cite{ABA06} Ben Abda et al introduce an energy functional, which
depends on both unknown traces of the solution of (CP) at $\Gamma_2$.

Our main goal in this article is to apply level set methods for
obtaining stable approximations of the solution of (CP).
The first application of level set methods to inverse problems was
proposed by Santosa \cite{San96}.
These methods can be used in identification problems where the unknown
parameter is a piecewise constant function assuming one of only two
possible values.
The methods considered here are adequate to solve (CP) in the special
case where the Neumann trace of $u$ at $\Gamma_2$ is known {\em a priori}
to satisfy $u_\nu |_{\Gamma_2} = \chi_D$, for some $D \subset \Gamma_2$.

A related application corresponds to the inverse problem in corrosion
detection \cite{CK07, Ing97}. This problem consists in determining information
about corrosion occurring on the inaccessible boundary part $(\Gamma_2$) of
a specimen. The data for this inverse problem correspond to prescribed
current flux ($g_2$) and voltage measurements ($g_1$) on the accessible
boundary part ($\Gamma_1$) and the model is the Laplace equation with no
source term ($P = \Delta$, $f = 0$). For simplicity one assumes the specimen
to be a thin plate ($\Omega \subset \R^2$) and $\partial\Omega = \Gamma_1
\cup \Gamma_2$. Moreover, the unknown corrosion damage $\gamma$ is assumed
to be the characteristic function of some $D \subset \Gamma_2$,
corresponding to the boundary condition $u_\nu + \gamma u = 0$ at
$\Gamma_2$.

The numerical methods analyzed in this article can be extended in a
straightforward way to arbitrary elliptic Cauchy problems possessing a
solution with similar structure, i.e. whenever the assumption that
$u_\nu |_{\Gamma_2}$ is a piecewise constant function assuming one of
only two possible values (not necessarily zero and one) is valid.
The proposed methods are inspired from the approaches followed in
\cite{FSL05,Bur01}, and relate to evolution flows of Hamilton-Jacobi type.
\medskip

The manuscript is outlined as follows: In Section~\ref{sec:2} we write the
elliptic Cauchy problem in the functional analytical framework of an
(ill-posed) operator equation. This is the starting point for the level set
approaches derived in the two subsequent sections.
In Section~\ref{sec:3} we investigate a level set method for (CP) based on
the approach proposed in \cite{Bur01}. We prove convergence and stability
of the proposed method, as well as a monotonicity result analog to the one
known for the asymptotic regularization method \cite{Tau94}.
In Section~\ref{sec:4} we derive a second level set method for (CP) based
on the ideas presented in \cite{FSL05}. First we prove existence of minima
for a least square functional related to (CP).
In the sequel we prove convergence and stability results for our
regularization strategy. The corresponding level set method is derived from
an explicit Euler method for solving the evolution equation related to
the first order optimality condition of the least square functional.
Section~\ref{sec:5} is devoted to numerics. Three different experiments
are provided, in order to illustrate the effectiveness of the level set
method considered in Section~\ref{sec:4}.

\section{Formulation of the inverse problem} \label{sec:2}

In this section we rewrite the elliptic Cauchy problem (CP) in the functional
analytical framework of an operator equation. This is the starting point for
the level set approaches derived in the Sections~\ref{sec:3}
and~\ref{sec:4}.

The functional analytical framework established in this section is similar to
the one derived in \cite{Le00}. The difference is that, instead of looking
for a fixed point operator, we follow an optimal control approach proposed
by Lions \cite{Lio71}. We begin by defining the auxiliary problem:
\begin{equation} \label{def:oper-T}
\left\{ \begin{array}{rl}
     P v   = f     \, ,&\!\!\! \mbox{in } \Omega \\
     v     = g_1   \, ,&\!\!\! \mbox{at } \Gamma_1 \\
     v_\nu = \vphi \, ,&\!\!\! \mbox{at } \Gamma_2
   \end{array} \right. .
\end{equation}
This mixed BVP defines the operator $T: \vphi \mapsto v_\nu |_{\Gamma_1}$.
Notice that, if $\vphi = u_\nu |_{\Gamma_2}$, where $u$ is the solution
of (CP), then it would follow $T(\vphi) = g_2$. Thus, a simple least square
approach \cite{CN06, Ing97} for (CP) consists in solving the optimization
problem.
$$  \| T(\vphi) - g_2 \|^2 \to \min . $$

Due to the superposition principle for linear elliptic BVPs \cite{GT77},
one can split the solution of (\ref{def:oper-T}) in $v = v_a + v_b$, where
%
%
\begin{equation} \label{def:va}
  P v_a     = 0     \, , \ \mbox{in } \Omega \quad\quad
  v_a       = 0     \, , \ \mbox{at } \Gamma_1 \quad\quad
  (v_a)_\nu = \vphi \, , \ \mbox{at } \Gamma_2 \, ;
\end{equation}
\begin{equation} \label{def:vb}
  P v_b     = f     \, , \ \mbox{in } \Omega \quad\quad
  v_b       = g_1   \, , \ \mbox{at } \Gamma_1 \quad\quad
  (v_b)_\nu = 0     \, , \ \mbox{at } \Gamma_2 \, .
\end{equation}
From (\ref{def:va}) we can define the linear operator
\begin{equation} \label{def:oper-L}
L: \vphi \mapsto (v_a)_\nu |_{\Gamma_1} \, ,
\end{equation}
and from (\ref{def:vb}) we define the function $z := (v_b)_\nu |_{\Gamma_1}$.
Since $T(\vphi) = L\, \vphi + z$, the Cauchy problem (CP) can be written
in the form of the operator equation
\begin{equation} \label{eq:op-cp}
L\, \vphi \ = \ g_2 - z \; ,
\end{equation}
where the constant term $z$ depends only on $g_1$, $f$ and $P$.
Therefore, it can be computed {\em a priori}.

In order to derive our first level set approach (Section~\ref{sec:3}) for
solving (\ref{eq:op-cp}), we have to formulate (CP) in such a way that $L$
is continuous with respect to the $L^2$-norm.
For the second level set approach (Section~\ref{sec:4}), we state (CP)
such that $L$ is continuous with respect to the $L^1$-norm.
In the sequel we present these two possible formulations of (CP).

\subsection{Framework for the first level set approach} \label{ssec:2-1}

We consider (CP) in the form of equation (\ref{eq:op-cp}). Further we 
assume the Cauchy data to satisfy
\begin{equation} \label{eq:op2-sobol-sp}
(g_1, g_2) \ \in \ H^{\me}(\Gamma_1) \times [H^{\me}_{00}(\Gamma_1)]'
\end{equation}
and the source term $f$ to be a $L^2(\Omega)$ distribution.

From this choice of $g_1$ and $f$, it follows that the mixed BVP in
(\ref{def:vb}) has a unique solution $v_b \in H^1(\Omega)$ \cite{GT77, Le00}.
Therefore, $z := (v_b)_\nu |_{\Gamma_1} \in [H^{\me}_{00}(\Gamma_1)]'$ and
the distribution $g_2 - z$ on the right hand side of (\ref{eq:op-cp}) is in
$[H^{\me}_{00}(\Gamma_1)]'$.

Notice that, if we choose $\vphi \in L^2(\Gamma_2) \subset
[H^{\me}_{00}(\Gamma_2)]'$, the mixed BVP in (\ref{def:va}) has a unique
solution $v_a \in H^1(\Omega)$ and the linear operator $L$ in
(\ref{def:oper-L}) is well defined from $L^2(\Gamma_2)$ into
$[H^{\me}_{00}(\Gamma_1)]'$. Indeed, this assertion follows from

\begin{propo} \label{prop:op2-L}
Let the domain $\Omega \subset \R^d$ with $d = 2, 3$ and the operator
$P$ be defined as in Section~\ref{sec:1}.
Then, the linear operator defined in (\ref{def:oper-L}) is an injective
bounded map $L : L^2(\Gamma_2) \to [H^{\me}_{00}(\Gamma_1)]'$.
\end{propo}
\proof
Since $d=2, 3$, the boundary part $\Gamma_2$ is either a 1D or a 2D Lipschitz
manifold and the embedding $L^2(\Gamma_2) \subset [H^{\me}_{00}(\Gamma_2)]'$
is continuous.
Therefore, given $\vphi \in L^2(\Gamma_2)$, the mixed BVP in (\ref{def:va})
has a unique solution $v_a \in H^1(\Omega)$ satisfying the {\em a priori}
estimate
$$ \| v_a \|_{H^1(\Omega)} \ \le \ C_1
   \| \vphi \|_{[H^{\me}_{00}(\Gamma_2)]'} \; , $$
for some positive constant $C_1$ (depending on $P$, $\Omega$ and $\Gamma_2$). 
Now, from the continuity of the Neumann trace operator $\gamma_{N,1}:
H^1(\Omega) \ni v \mapsto v_\nu|_{\Gamma_1} \in [H^{\me}_{00}(\Gamma_1)]'$,
follows
$$ \| L \vphi \|_{[H^{\me}_{00}(\Gamma_1)]'} \ \le \ C_2
   \| v_a \|_{H^1(\Omega)} \ \le \ C_3
   \| \vphi \|_{L^2(\Gamma_2)} \; , $$
and the continuity of $L$ follows. It remains to prove the injectivity
of $L$.
Notice that, if $L \vphi = 0$, the function $v_a$ in (\ref{def:va})
satisfies: $P v_a = 0$ in $\Omega$, $v_a = (v_a)_\nu = 0$ at $\Gamma_1$.
Then, $\vphi = 0$ follows from the uniqueness of weak solution for (CP)
\cite{EL01}. \hfill \endproof

Summarizing, this setup allow us to state (CP) in the form of the operator
equation (\ref{eq:op-cp}), where $L$ is the linear continuous operator
\begin{equation} \label{eig:op1-L}
L : L^2(\Gamma_2) \to [H^{\me}_{00}(\Gamma_1)]'
\end{equation}
defined in (\ref{def:oper-L}).

\subsection{Framework for the second level set approach} \label{ssec:2-2}

In the sequel we shall assume $\Omega \subset \R^d$ with $d = 2,3$ and define
yet another functional analytical framework to analyze (\ref{eq:op-cp}).
The Cauchy data is assumed to satisfy
\begin{equation} \label{eq:op1-sobol-sp}
(g_1, g_2) \ \in \ [H^{\me}_{00}(\Gamma_1)]' \times
                   [H^{\tme}_{00}(\Gamma_1)]'
\end{equation}
and the source term $f$ to be a $H^{-1}(\Omega)$ distribution.

Due to the choice of $g_1$ and $f$ above, the elliptic theory allow us
to conclude that the mixed BVP in (\ref{def:vb}) has a unique solution
$v_b \in L^2(\Omega)$ \cite{DL88, GT77}.
Therefore, $z := (v_b)_\nu |_{\Gamma_1} \in [H^{\tme}_{00}(\Gamma_1)]'$
and the term $g_2 - z$ on the right hand side of (\ref{eq:op-cp}) is a
distribution in $[H^{\tme}_{00}(\Gamma_1)]'$.

In the next proposition we prove that the linear operator $L$ in
(\ref{def:oper-L}) is well defined continuous and injective from
$L^1(\Gamma_2)$ to $[H^{\tme}_{00}(\Gamma_1)]'$.

\begin{propo} \label{prop:op1-L}
Let the domain $\Omega \subset \R^d$ with $d = 2,3$ and the operator $P$
be defined as in Section~\ref{sec:1}.
Then, the linear operator defined in (\ref{def:oper-L}) is an injective
bounded map $L : L^1(\Gamma_2) \to [H^{\tme}_{00}(\Gamma_1)]'$.
\end{propo}
\proof
If $d=2$, the boundary part $\Gamma_2$ is a 1D manifold and the Sobolev
embedding theorem \cite{Ada75, GT77} implies
$H^s_0(\Gamma_2) \subset L^\infty(\Gamma_2)$ for $s > 1/2$.
If $d=3$, the embedding above still holds, however only for $s > 1$.
In either case, we can take $s = \tme$ and conclude that
$L^1(\Gamma_2) \subset [L^\infty(\Gamma_2)]' \subset H^{-\tme}(\Gamma_2)$.
This result together with
$H^{-\tme}(\Gamma_2) \subset [H^{\tme}_{00}(\Gamma_2)]'$ imply the
continuity of the embedding $L^1(\Gamma_2) \subset [H^{\tme}_{00}(\Gamma_2)]'$.
Thus, from the elliptic theory we conclude that given $\vphi \in
L^1(\Gamma_2)$, the mixed BVP in (\ref{def:va}) has a unique solution
$v_a \in L^2(\Omega)$ satisfying the {\em a priori} estimate
$$ \| v_a \|_{L^2(\Omega)} \ \le \ C_1
   \| \vphi \|_{[H^{\tme}_{00}(\Gamma_2)]'} \; , $$
for some positive constant $C_1$ (depending on $P$, $\Omega$ and $\Gamma_2$). 
Now, from the continuity of the Neumann trace operator $\gamma_{N,1}:
L^2(\Omega) \ni v \mapsto v_\nu|_{\Gamma_1} \in [H^{\tme}_{00}(\Gamma_1)]'$,
follows
$$ \| L \vphi \|_{[H^{\tme}_{00}(\Gamma_1)]'} \ \le \ C_2
   \| v_a \|_{L^2(\Omega)} \ \le \ C_3
   \| \vphi \|_{L^1(\Gamma_2)} \; , $$
and the continuity of $L$ is proven. In order to prove the injectivity of
$L$, one argues analogously as in the last part of the proof of Proposition~%
\ref{prop:op2-L}.
\hfill \endproof

Summarizing, this setup allow us to state (CP) in the form of equation
(\ref{eq:op-cp}), where $L$ is the linear continuous operator
\begin{equation} \label{eig:op2-L}
L : L^1(\Gamma_2) \to [H^{\tme}_{00}(\Gamma_1)]'
\end{equation}
defined in (\ref{def:oper-L}).

\subsection{A remark on the dimension of $\Omega \subset \R^d$}
\label{ssec:2-3}

In the approach presented in Subsection~\ref{ssec:2-1}, the key argument to
prove the desired regularity of the operator $L$ (see (\ref{eig:op1-L})) was
the continuity of the embedding $L^2(\Gamma_2) \subset
[H^{\me}_{00}(\Gamma_2)]'$.
Since the continuity of this embedding does not depend on $d$,
Proposition~\ref{prop:op2-L} actually holds for every $d \in \N$.

In the case of the approach presented in Subsection~\ref{ssec:2-2}, the
situation is different:
Proposition~\ref{prop:op1-L} only holds for $d = 2, 3$. Indeed, the desired
regularity of $L$ in (\ref{eig:op2-L}) depends on the continuity of the
embedding $L^1(\Gamma_2) \subset [H^{\tme}_{00}(\Gamma_2)]'$, which
follows from the Sobolev embedding theorem:
$$ H^s_0(\Gamma_2) \subset C^0(\overline{\Gamma}_2) \subset
   L^\infty(\Gamma_2) \, , \ {\rm for } \ s > (d-1)/2 $$
(notice that $\Gamma_2 \in \R^{d-1}$).
Therefore, if $\Omega \in \R^d$ for $d \ge 4$, we would have to choose
$s \ge 2$ in the proof of Proposition~\ref{prop:op1-L}. The proof would
no longer hold, since the mixed BVP in (\ref{def:va}) would not have a
solution $v_a$ on a Sobolev space of non-negative index.

This is not actually a disadvantage of the second level set approach,
since almost all real life applications are related to either two or
three-dimensional domains $\Omega$.

\subsection{A remark on noisy Cauchy data} \label{ssec:2-4}

The second remark concerns the investigation of (CP) for noisy data.
If only corrupted noisy data $(g_1^\delta, g_2^\delta)$ are available,
we assume the existence of a consistent pair of Cauchy data $(g_1, g_2)$
such that
\begin{equation} \label{eq:assump-noise}
\| g_1 - g_1^\delta \|_{L^2(\Gamma_1)} +
\| g_2 - g_2^\delta \|_{L^2(\Gamma_1)} \ \le \ \delta \, .
\end{equation}

For clarity of the presentation, we discuss separately the two frameworks
introduced above in this section:
\medskip

\noindent
{\bf Noisy data within the framework of Subsection~\ref{ssec:2-2}:} \\
Let the noisy data be given as in (\ref{eq:assump-noise}) and the exact
Cauchy data satisfy (\ref{eq:op1-sobol-sp}). Since $z$ in (\ref{eq:op-cp})
depends continuously on $g_1$ in the $[H^{\me}_{00}(\Gamma_1)]'$ topology,
we can solve the mixed BVP in (\ref{def:vb}) using $g_1^\delta$ as data, and
obtain a corresponding $z^\delta \in [H^{\tme}_{00}(\Gamma_1)]'$ such that
\begin{equation} \label{eq:noise-rewritten}
\| (g_2 - z) - (g_2^\delta - z^\delta) \|_{[H^{\tme}_{00}(\Gamma_1)]'}
   \ \le \ C \delta \, ,
\end{equation}
where the constant $C$ depends on $\Omega$, $P$, $\Gamma_1$ and $f$.
Summarizing, we have

\begin{lemma} \label{lem:sumary2}
Within the framework of Subsection~\ref{ssec:2-2}, the Cauchy problem
(CP) with noisy data satisfying (\ref{eq:assump-noise}) reduces to
equation
$$ L\, \vphi \ = \ g_2^\delta - z^\delta \; , $$
where $L$ satisfies (\ref{eig:op2-L}) and the right hand side
$(g_2^\delta - z^\delta)$ satisfies (\ref{eq:noise-rewritten}).
\end{lemma}

\noindent
{\bf Noisy data within the framework of Subsection~\ref{ssec:2-1}:} \\
Let the noisy data be given as in (\ref{eq:assump-noise}) and the exact
Cauchy data satisfy (\ref{eq:op2-sobol-sp}). 
Since $z$ in (\ref{eq:op-cp}) depends continuously on $g_1$ in the
$H^{\me}(\Gamma_1)$ topology, a natural question arises:

\begin{minipage}{12cm}
{\it Q) Is it possible to obtain from measured data $(g_1^\delta, g_2^\delta)$
satisfying (\ref{eq:assump-noise}), a corresponding $z^\delta \in
[H^{\me}_{00}(\Gamma_1)]'$ such that
$\| z - z^\delta \|_{[H^{\me}_{00}(\Gamma_1)]'}  \le  \delta$?}
\end{minipage}

\noindent
We claim that such $z^\delta$ can be obtained under the {\em a priori}
assumption $g_1 \in H^{s}(\Gamma_1)$, for some $s > \me$. Indeed, under this
assumption, \cite[Lemma~8]{EL01} guarantees the existence of a smoothing
operator $S: L^2(\Gamma_1) \to H^{\me}(\Gamma_1)$, and of a positive function
$\mu$ with $\lim_{t \downarrow 0} \mu(t) = 0$, such that for $\delta > 0$
and $g_1^\delta \in L^2(\Gamma_1)$ with
$\| g_1 - g_1^\delta \|_{L^2(\Gamma_1)} \le \delta$, we have
$\| g_1 - S(g_1^\delta) \|_{H^{\me}(\Gamma_1)} \le \mu(\delta)$.
Thus, after smoothing the noisy data $g_1^\delta$, we obtain
$\hat{g}_1^\delta := S(g_1^\delta) \in H^{\me}(\Gamma_1)$. Next, we
solve the mixed BVP in (\ref{def:vb}) using $\hat{g}_1^\delta$ as data,
and obtain a corresponding $z^\delta \in [H^{\me}_{00}(\Gamma_1)]'$ with
$\| z - z^\delta \|_{[H^{\me}_{00}(\Gamma_1)]'}  \le  C \mu(\delta)$,
with $C$ as in (\ref{eq:noise-rewritten}). \\
Once we are able to give an affirmative answer to question Q), it
follows from (\ref{eq:assump-noise}) that
\begin{equation} \label{eq:noise-rewritten2}
\| (g_2 - z) - (g_2^\delta - z^\delta) \|_{[H^{\me}_{00}(\Gamma_1)]'}
   \ \le \ (1+C) \max\{ \delta, \mu(\delta) \} \, .
\end{equation}
Summarizing, we have

\begin{lemma} \label{lem:sumary1}
Consider the framework of Subsection~\ref{ssec:2-1}. Assume the noisy
Cauchy data satisfy (\ref{eq:assump-noise}), where $g_1 \in H^s(\Gamma_1)$
for some $s > 1/2$. Then (CP) reduces to equation
$$ L\, \vphi \ = \ g_2^\delta - z^\delta \; , $$
where $L$ satisfies (\ref{eig:op1-L}) and the right hand side
$(g_2^\delta - z^\delta)$ satisfies (\ref{eq:noise-rewritten2}).
\end{lemma}

\section{A first level set approach} \label{sec:3}

In this section we investigate a level set method for (CP) based on the
approach introduced in \cite{Bur01}. For this purpose, we consider the
functional analytical framework derived in Subsection~\ref{ssec:2-1}
and summarized in Lemma~\ref{lem:sumary1}.

The starting point of this approach is the assumption that the solution
$\overline{\vphi}$ of (\ref{eq:op-cp}) is the characteristic function
$\chi_{D}$ of a subdomain $D \subset\subset \Gamma_2$.
A function $\phi: \Gamma_2 \times \R^+ \to \R$ is introduced, allowing the
definition of the level sets $D(t) = \{ \phi(\cdot,t) \geq 0 \}$. The
function $\phi$ should be chosen such that
\begin{equation} \label{def:vphi}
\vphi(\cdot,t) := \chi_{D(t)} \to \chi_{D} = \overline{\vphi}
\end{equation}
as $t \to \infty$. The level set method corresponds to a continuous evolution
for an artificial time $t$, where the {\em level set function} $\phi$ is
defined by a Hamilton-Jacob equation of the form
\begin{equation} \label{eq:HJ}
\frac{\partial \phi}{\partial t} + V \cdot \nabla \phi \ = \ 0 \, ,
\end{equation}
with initial value $\phi(x,0) = \phi_0(x)$, where $\phi_0$ is an appropriate
indicator function of a measurable set $D_0 \subset\subset \Gamma_2$.
The function $V: \R^{d-1} \times \R^+ \to \R^{d-1}$ describes the velocity
of the level sets of $\phi$. Following \cite{San96}, $V$ is chosen in the
normal direction of the level set curves of $\phi$, i.e. $V =
v \, \nabla\phi/|\nabla\phi|$, for $v: \R^{d-1} \times \R^+ \to \R$.

Here, the guideline for the choice of the {\em velocity function} $V$ is
a property of the asymptotic regularization method \cite{Tau94}. In this
method the approximations $\vphi(\cdot,t)$ for the solution
$\overline{\vphi}$ of (\ref{eq:op-cp}) satisfy
\begin{equation} \label{eq:assymp-reg}
\frac{d}{d t} \| \varphi(\cdot,t) - \overline{\vphi} \|^2  =
- 2 \| L (\varphi(\cdot,t)) - (g_2 - z) \|^2 \, .
\end{equation}
As we shall see, a level set method satisfying \eqref{eq:assymp-reg} can
be analyzed in a similar way to asymptotic regularization.%
\footnote{The method of asymptotic regularization can not be used directly
to construct piecewise approximations for $\overline{\varphi}$, since it
uses the time derivative of $\varphi(\cdot, t)$ which is not defined in
$L^2(\Gamma_2)$.}
Our first goal is to determine how to choose the velocity $V$ such that,
if $\phi$ solves (\ref{eq:HJ}) then the function $\vphi$ defined in
(\ref{def:vphi}) satisfies (\ref{eq:assymp-reg}).

\begin{propo} \label{prop:deriv-bur}
Let the function $h$ be defined by
$$ h(x, t) = (-1+2\varphi(x, t)) \cdot {\rm div}\, V(x, t) , \
   x \in \Gamma_2 \, , \ t \in \R^+ , $$
for some $V \in L^{\infty} (0, T, L^2(\R^{d-1}))^{d-1}$ with
${\rm div}\, V \in L^1 (0, T, L^{\infty}(\R^{d-1}))$  $\cap$
$L^{\infty} (0, T, L^2(\R^{d-1}))$. Moreover, let the level set function
$\phi$ satisfy \eqref{eq:HJ} and $\vphi$ be defined by (\ref{def:vphi}).
Then
\begin{equation} \label{eq:deriv-bur}
\frac{d}{d t} \| \vphi(\cdot,t) - \overline{\varphi} \|_{L^2(\Gamma_2)}^2
\ = \
\int_{\Gamma_2} (\vphi(\cdot, t) - \overline{\vphi}) \cdot h(x,t) dx\, ,
\end{equation}
for all $t \in (0,T)$.
\end{propo}
\proof See \cite[Proposition~3.3]{Bur01}.
\hfill \endproof

\begin{remark}
If $V$ satisfies the assumptions of Proposition~\ref{prop:deriv-bur} for
each $T > 0$, then (\ref{eq:deriv-bur}) is satisfied for all $t \in \R^+$.
\end{remark}

Since we are assuming the operator $L$ to satisfy (\ref{eig:op1-L}), we
conclude from Proposition~\ref{prop:deriv-bur} that relation
(\ref{eq:assymp-reg}) is satisfied for the velocity $V$ satisfying
\begin{equation} \label{eq:choice-V}
- {\rm div}\, V =
2 (-1 + 2 \vphi)^{-1} \, L^* ( L\, \vphi - (g_2 - z) ) \ 
{\rm in } \ \Gamma_2 \times \R^{+} \, .
\end{equation}
Observe that, if the normal derivative $V_\nu(\cdot,t)$ vanishes on
$\partial\Gamma_2$ for $t \in \R^+$ then the support of $\vphi(\cdot,t)$
remains a subset of $\overline{\Omega}$ during the evolution, what is a
desirable property.
It is worth noticing that a solution $V$ of (\ref{eq:choice-V}) with
homogeneous Neumann boundary condition on $\partial\Gamma_2$ always exists.
Indeed, it is enough to choose $V = \nabla\psi$, where $\psi$ solves
$$ -\Delta \psi \ = \ 2 (-1 + 2 \vphi)^{-1} \, L^* (L\, \vphi - (g_2 - z))\, ,
                    \ {\rm on } \ \Gamma_2\, , \quad
   \psi \ = \ 0\, , \ {\rm at } \ \partial\Gamma_2\, . $$

In the sequel we derive a convergence analysis for the level set method
defined by (\ref{def:vphi}) (\ref{eq:HJ}), with the choice of velocity
in (\ref{eq:choice-V}). We consider noisy Cauchy data as in Lemma~%
\ref{lem:sumary1}. Moreover, we define the {\em stopping time}
$T(\delta, g_1^\delta, g_2^\delta)$ by the generalized discrepancy
principle \cite{EHN96}
\begin{equation} \label{eq:discrep}
T(\delta, g_1^\delta, g_2^\delta) \ := \
\inf \{ t \in \R^+; \ \| L(\vphi(\cdot, t)) - (g_2^\delta - z^\delta) \|
        \leq \tau \delta \} \, ,
\end{equation}
for some $\tau > 1$.

The next theorem summarizes the main convergence and stability results
for this level set method

\begin{theorem}[Convergence analysis] \label{th:conv-anal}
Let $V$ satisfy (\ref{eq:choice-V}) and $\vphi$, $\phi$ be defined by
(\ref{def:vphi}) (\ref{eq:HJ}).
\begin{itemize}
\item[i)] {\bf Monotonicity:}
For noisy Cauchy data and $\tau > 1$, the iteration error is strictly
monotone decreasing, i.e.
$$ \frac{d}{dt} \| \vphi(\cdot, t) - \overline{\vphi}\|^2 < 0 \, , $$
for all $t > 0$ with
$\| L (\vphi(\cdot, t)) - (g_2^\delta - z^\delta) \| > \tau \delta$.
Moreover, for exact Cauchy data (i.e., $\delta=0$), we have the inequality
$$ \int_0^\infty \| L(\vphi(\cdot, t)) - L(\overline{\vphi})\|^2 \, dt
   < \infty \, ; $$
\item[ii)]  {\bf Convergence:}
If the Cauchy data is exact, then $\vphi(\cdot, t) \to \overline{\vphi}$
in $L^2(\Gamma_2)$ as $t \to \infty$, where $\overline{\vphi}$ is the
solution of (\ref{eq:op-cp}) corresponding to the consistent data $(g_1,g_2)$;
\item[iii)]  {\bf Stability:}
For noisy Cauchy data, the stopping time $T(\delta,g_1^\delta,g_2^\delta)
=: T_\delta$ defined by (\ref{eq:discrep}) with $\tau > 1$ is finite.
Moreover, given a sequence $\delta_k \to 0$ and
$\{ (g_1^{\delta_k}, g_2^{\delta_k}) \}_k$ corresponding noisy data satisfying
(\ref{eq:assump-noise}) for some consistent data pair $(g_1,g_2)$, then the
approximations $\vphi(\cdot, T_{\delta_k})$ converge to $\overline{\vphi}$
in $L^2(\Gamma_2)$ as $\delta_k \to 0$, where $\overline{\vphi}$ is the
solution of (\ref{eq:op-cp}) corresponding to $(g_1,g_2)$.
\end{itemize}
\end{theorem}
\proof
Item {\it i)}: The monotonicity result is a consequence of
(\ref{eq:assymp-reg}). The inequality in the second statement is a well
known property of the asymptotic regularization and its proof is analog
to \cite[Proposition~4.1]{Bur01}. \\
Item {\it ii)}: This statement concerning convergence for exact data
is also known to hold for the Landweber iteration as well as for asymptotic
regularization. The proof follow the lines of \cite[Theorem~4.4]{Bur01}. \\
Item {\it iii)}: This convergence result for noisy data has a counterpart
in the asymptotic regularization. The proof carries over from
\cite[Theorem~4]{Tau94}.
\hfill\endproof

\begin{remark}
Under the assumptions of Theorem~\ref{th:conv-anal} one can prove,
analogous as in item i), that the residual function
$t \mapsto \| L(\vphi(\cdot, t)) - (g_2^\delta - z^\delta) \|$ is
monotonically non-increasing.
\end{remark}

\section{A second level set approach} \label{sec:4}

In the sequel we investigate a second level set method for (CP). Our
approach is based on \cite{FSL05}.
In what follows we shall consider the functional analytical framework
for (CP) discussed in Subsection~\ref{ssec:2-2} and summarized in
Lemma~\ref{lem:sumary2}.

Let functions $\vphi$ and $\phi$ be defined as in Section~\ref{sec:3}.
For simplicity we adopt the notation $Y := H^{3/2}_{00}(\Gamma_1)'$.
If we denote by $H$ the heavy-side projector%
\footnote{The projector $H$ and its approximation $H_\eps$ are defined by
$$ H(t) := \left\{ \begin{array}{rcl}
             0 & \mbox{ for } & t < 0 \, , \\
             1 & \mbox{ for } & t \geq 0 \; .
           \end{array} \right. \quad \quad
H_\eps(t) := \left\{ \begin{array}{rcl}
                0          & \mbox{ for } & t < -\eps \, , \\
                1 + t/\eps & \mbox{ for } & -\eps \le t \le 0 \; . \\
                1          & \mbox{ for } & t \geq 0 \; .
             \end{array} \right. $$ }
then the Cauchy problem (\ref{eq:op-cp}) can be written in the form of the
constrained optimization Problem
$$ \min \| L\, \vphi - (g_2^\delta - z^\delta) \|^2_Y \, , \quad
   {\rm s.t. } \quad \vphi = H(\phi) \, . $$
Alternatively, we can minimize
\begin{equation} \label{eq:cp-lset}
 \min \| L (H(\phi)) - (g_2^\delta - z^\delta) \|^2_Y \, ,
\end{equation}
over $\phi \in H^1(\Gamma_2)$. Tikhonov regularization for (\ref{eq:cp-lset})
using $TV-H^1$ penalization consists in the minimization of the cost
functional
\begin{equation} \label{eq:lset-tik}
{\cal F}_\alpha (\phi) := \| L(H(\phi)) - (g_2^\delta - z^\delta) \|^2_Y +
\alpha \big[ \beta | H(\phi) |_{BV} + \|\phi - \phi_0\|_{H^1}^2 \big] \, ,
\end{equation}
where $\alpha > 0$ plays the rule of a regularization parameter and
$\beta > 0$ is a scaling factor \cite{LS03, FSL05}.
Since $H$ is a discontinuous operator, one cannot prove that the
Tikhonov functional in (\ref{eq:lset-tik}) attains a minimizer.

In order to guarantee existence of a minimizer of ${\cal F}_\alpha$,
we use the concept of generalized minimizers in \cite[Lemma~2.2]{FSL05}.
${\cal F}_\alpha$ is no longer considered as a functional on $H^1$, but
as a functional defined on the w-closure of the graph of $H$, contained
in $BV \times H^1(\Gamma_2)$.
A generalized minimizer of ${\cal F}_\alpha(\phi)$ is defined as a
minimizer of
\begin{equation} \label{eq:lset-tik-fls}
{\cal F}_\alpha (\xi,\phi) := \| L(H(\phi)) - (g_2^\delta - z^\delta) \|^2_Y
                              + \alpha \, \rho(\xi,\phi)
\end{equation}
on the set of admissible pairs
\begin{multline*}
Ad \, := \, \{ (\xi,\phi) \in L^\infty(\Gamma_2) \times H^1(\Gamma_2) \, ; \
\exists \ \{\phi_k\} \in H^1 \ {\rm and } \ \{ \eps_k \} \in {\mathbb R}^+ \
{\rm s.t. } \\
\lim_{k\to\infty} \eps_k = 0 \, , \
\lim_{k\to\infty} \| \phi_k - \phi \|_{L^2} = 0 \, , \
\lim_{k\to\infty} \| H_{\eps_k}(\phi_k) - \xi \|_{L^1} = 0 \} \, ,
\end{multline*}
where $\displaystyle \rho(\xi,\phi) := \inf_{\{\phi_k\}, \{\eps_k\}}
\liminf_{k \to \infty} \left\{ 2 \beta | H_{\eps_k} (\phi_k)|_{BV} +
\|\phi_k-\phi_0\|_{H^1}^2 \right\}$.

As a consequence of this definition, the penalization term in
(\ref{eq:lset-tik}) can be interpreted as a functional $\rho:
Ad \to \mathbb R^+$. In order to prove coerciveness and weak lower
semi-continuity of $\rho$, the assumption that $L$ is a continuous
operator on a $L^1$ space is crucial (see Proposition~\ref{prop:op1-L}).
These properties of $\rho$ are the main arguments needed to prove
existence of a generalized minimizer $(\overline{\xi}_\alpha,
\overline{\phi}_\alpha)$ of ${\cal F}_\alpha$ in $Ad$
\cite[Theorem~2.9]{FSL05}.

The classical analysis of Tikhonov type regularization methods \cite{EHN96}
do apply to functional ${\cal F}_\alpha$, as we shall see next.

\begin{theorem}[Convergence] \label{th:conv-anal2}
Let $\alpha: \mathbb R^+ \to \mathbb R^+$ be a function satisfying
$\displaystyle\lim_{\delta \to 0} \alpha (\delta) = 0$ and
$\displaystyle\lim_{\delta \to 0} \delta^2 \, \alpha^{-1}(\delta) = 0$.
Given a sequence $\delta_k \to 0$ and
$\{ (g_1^{\delta_k}, g_2^{\delta_k}) \}_k$ corresponding noisy data
satisfying (\ref{eq:assump-noise}) for some consistent data pair
$(g_1,g_2)$, then the minimizers $(\xi_k,\phi_k)$ of
${\cal F}_{\alpha(\delta_k)}$ converge in $L^1(\Gamma_2) \times
L^2(\Gamma_2)$ to a minimizer $(\overline{\xi}_\alpha,
\overline{\phi}_\alpha)$ of ${\cal F}_\alpha$ in (\ref{eq:lset-tik-fls}).
\end{theorem}
\proof
See \cite[Section~2.3]{FSL05}.
\hfill\endproof

\subsection*{Numerical realization}

What concerns numerical approximations to the solution of (\ref{eq:op-cp}),
the functional ${\cal F}_\alpha$ in (\ref{eq:lset-tik}) has an interesting
property. Namely, it's generalized minimizers can be approximated by
minimizers of the stabilized functional%
\addtocounter{footnote}{-1}\footnotemark\addtocounter{footnote}{1}
\begin{equation} \label{eq:lset-tik-eps}
{\cal F}_{\alpha,\eps} (\phi) := \| L(H_\eps(\phi)) - (g_2 - z) \|^2_Y +
\alpha \big[ \beta | H_\eps(\phi) |_{BV} + \|\phi - \phi_0\|_{H^1}^2 \big]\,.
\end{equation}
In other words: Let $\phi_{\alpha,\eps}$ be a minimizer of
${\cal F}_{\alpha,\eps}$; given a sequence $\eps_k \to 0^+$, we can find a
subsequence $( H(\phi_{\alpha,\eps}), \phi_{\alpha,\eps} )$ converging in
$L^1(\Gamma_2) \times L^2(\Gamma_2)$ and the limit minimizes
${\cal F}_{\alpha,}$ in $Ad$.

The existence of minimizers of ${\cal F}_{\alpha,\eps}$ in $H^1(\Gamma_2)$
still has to be cleared:
Since $H_\eps$ is continuous and the operator $L$ is linear, continuous,
satisfying (\ref{eig:op2-L}), the existence of minimizers for
${\cal F}_{\alpha,\eps}$ follows directly from \cite[Lemma~3.1]{FSL05}.

This relation between the minimizers of ${\cal F}_\alpha$ and
${\cal F}_{\alpha,\eps}$ is the starting point for the derivation of a
numerical method. We can formally write the optimality condition for
${\cal F}_{\alpha,\eps}$ as
$$ \alpha (I-\Delta)(\phi - \phi_0) = {\cal R}_{\alpha,\eps}(\phi) \; , $$
where
\begin{multline*}
{\cal R}_{\alpha,\eps}(\vphi) := 
 - H'_\eps(\phi) L'(H_\eps(\phi))^* [L(H_\eps(\phi)) - (g_2^\delta - z^\delta)]
 \\ + \beta \alpha H'_\eps(\phi) \nabla \cdot
   ( \nabla H_\eps(\phi) / |\nabla H_\eps(\phi)| ) .
\end{multline*}
Identifying $\alpha = 1/\Delta t$, $\phi(0) = \phi_0$, $\phi(\Delta t) =
\phi$, we find
$$ (I-\Delta) \left( \frac{\phi(\Delta t) - \phi(0)}{\Delta t} \right)
   = {\cal R}_{1/\Delta t,\eps}(\phi(\Delta t)) \; . $$
Considering $\Delta t$ as a time discretization, we find that (in a formal
sense) the iterative regularized solution $\phi(\Delta t)$ is a solution of
an implicit time step for the dynamic system
\begin{equation} \label{eq:evol}
(I-\Delta) \left( \frac{\partial \phi(t)}{\partial t} \right) =
R_{1/\Delta t,\eps} (\phi(t)) \; .
\end{equation}
Our second level set method is based on the solution of the dynamic system
(\ref{eq:evol}). In algorithmic form we have:

\begin{tt}\begin{itemize}
\item[{\bf 0.}] \ Choose $\phi_0 \in H^1(\Gamma_2)$ and set $k = 0$; Compute
                $z^\delta = u_\nu|_{\Gamma_1}$, where
$$ \Delta u = f , \ {\tt in}\ \Omega \, , \quad
   u|_{\Gamma_1} = g_1^\delta\, , \quad
   u_\nu|_{\Gamma_2} = 0\, ; $$
\item[{\bf 1.}] \ Evaluate the residual
                $r_k := L(H_\eps(\phi_k)) - (g_2^\delta - z^\delta)$.
                Notice that $L(H_\eps(\phi_k)) = (u_k)_\nu |_{\Gamma_1}$, where
$$ \Delta u_k = 0 , \ {\tt in}\ \Omega \, , \quad
   u_k|_{\Gamma_1} = 0\, , \quad
   (u_k)_\nu|_{\Gamma_2} = H_\eps(\phi_k)\, ; $$
\item[{\bf 2.}] \ Evaluate $V_k := L'( H_\eps(\phi_k) )^*(r_k)$. Notice that
                $(L')^*(r_k) = -v_k |_{\Gamma_2}$, where
$$ \Delta v_k = 0\, ,\ {\tt in}\ \Omega\, ; \quad
   v_k|_{\Gamma_1} = r_k\, , \quad
   (v_k)_\nu|_{\Gamma_2} = 0\, ; $$
\item[{\bf 3.}] \ Evaluate the velocity $w_k$ by solving
\[\begin{aligned}
(I - \Delta) w_k &= H'_\eps(\phi_k) \, \Big[ - v_k + \beta\, \nabla \cdot
   \Big( \frac{\nabla H_\eps(\phi_k)}{|\nabla H_\eps (\phi_k)|} \Big) \Big]
   \, , {\tt in}\ \Gamma_2 ; \\
   (w_k)_\nu |_{\partial\Gamma_2} & = 0 \; .
  \end{aligned}
\]
\item[{\bf 4.}] \ Update the level set function\, $\phi_{k+1} =
\phi_k + \frac{1}{\alpha} \; w_k$;
\end{itemize} \end{tt}

Notice that Steps~1 and~2 above involve the solution of a mixed BVP in
$\Omega \subset \mathbb R^d$. On the third step, the computation of the
velocity function for the level set method requires a solution of a
Neumann BVP at $\Gamma_2$. 

In the next section we present some numerical experiments, which where
implemented using the above algorithm.

\section{Numerical experiments} \label{sec:5}

In what follows we present three numerical experiments for the level set
method analyzed in Section~\ref{sec:4}.
In Subsection~\ref{ssec:5-1} we consider a problem with exact Cauchy data,
where the solution is the characteristic function of a non-connected set.
In Subsection~\ref{ssec:5-2} we investigate how the degree of ill-posedness
of an elliptic Cauchy problem affects the performance of the level set
method.
In Subsection~\ref{ssec:5-3} we consider a problem with noisy Cauchy data
and test the stability of our method.

\subsection{A Cauchy problem with non-connected solution} \label{ssec:5-1}

\begin{figure}[t]
\centerline{ \epsfxsize4.6cm \epsfysize4.0cm
             \epsfbox{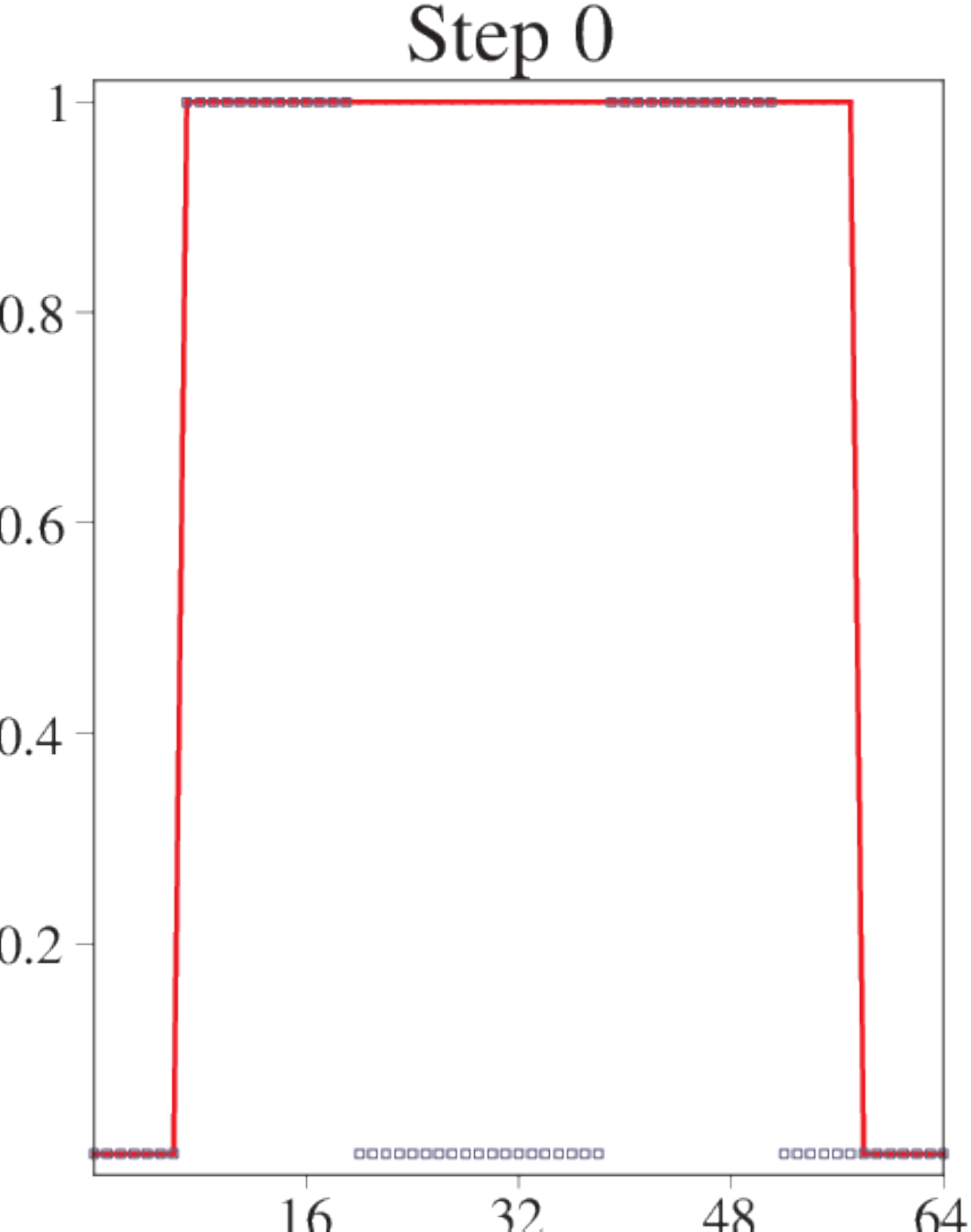} \hfil
             \epsfxsize6.3cm \epsfbox{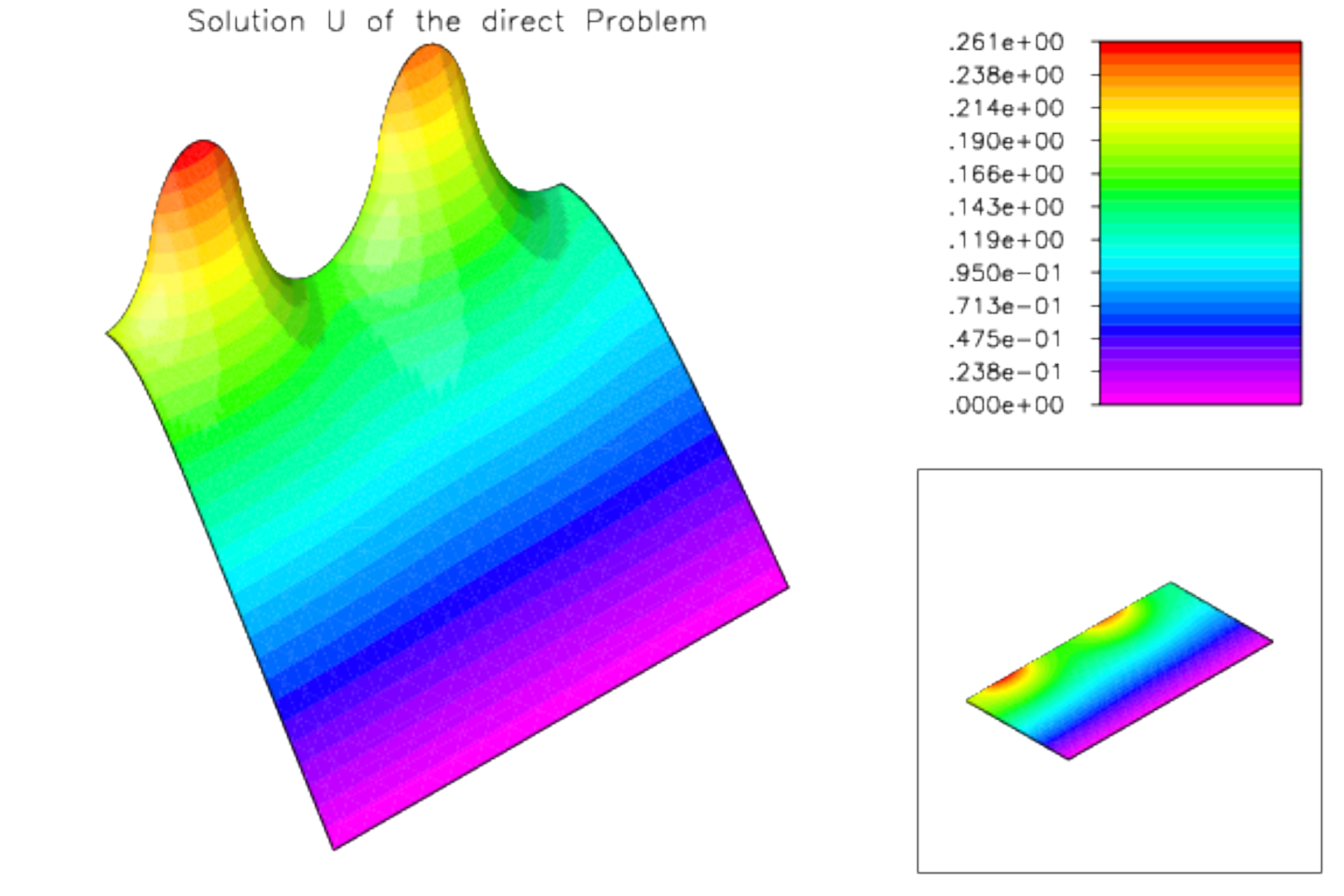}}
\caption{\footnotesize Framework for the first numerical experiment:
On the left hand side, the exact solution $\overline{\vphi}$ of the Cauchy
problem (dotted blue line) and the initial guess for the level set method
(solid red line).
The solution $u$ of the elliptic BVP corresponding to (\ref{eq:cp-exp1})
is depicted on the right hand side.
Notice that $u_\nu = g_2$ and $u = 0$ at $\Gamma_1$ (lower edge). Moreover,
$u_\nu = \overline{\vphi}$ at $\Gamma_2$ (top edge). \label{fig:exp1-frwork}}
\end{figure}

One of the main advantages of the level set approach is the fact that no
{\em a priori} assumption on the topology of the solution set is needed
\cite{San96}.
In this first numerical experiment we exploit this feature of the method
by solving a Cauchy problem with non-connected solution.
Let $\Omega = (0,1) \times (0,0.5)$ with $\Gamma_1 := (0,1) \times \{0\}$,
$\Gamma_2 := (0,1) \times \{0.5\}$ and $\Gamma_3 := \partial\Omega /
\{ \Gamma_1 \cup \Gamma_2\}$. Consider the Cauchy problem
\begin{equation} \label{eq:cp-exp1}
\Delta u  =  0 \, , \ \mbox{in } \, \Omega   \quad
 u  =  0       \, , \ \mbox{at } \, \Gamma_1 \quad
 u_\nu  =  g_2 \, , \ \mbox{at } \, \Gamma_1 \quad
 u_\nu  =  0   \, , \ \mbox{at } \, \Gamma_3 \, .
\end{equation}
The Cauchy data $(0,g_2)$ is chosen in such a way that the solution
$\overline{\vphi}$ of (\ref{eq:cp-exp1}) is the indicator function of
a non-connected subset of $\Gamma_2$ (see Figure~\ref{fig:exp1-frwork}).

\begin{figure}[t]
\centerline{ \epsfxsize4.2cm \epsfysize3.6cm
             \epsfbox{ps/seg2-sqr4-err0-ls-s00000} \hfil
             \epsfxsize4.2cm \epsfysize3.6cm
             \epsfbox{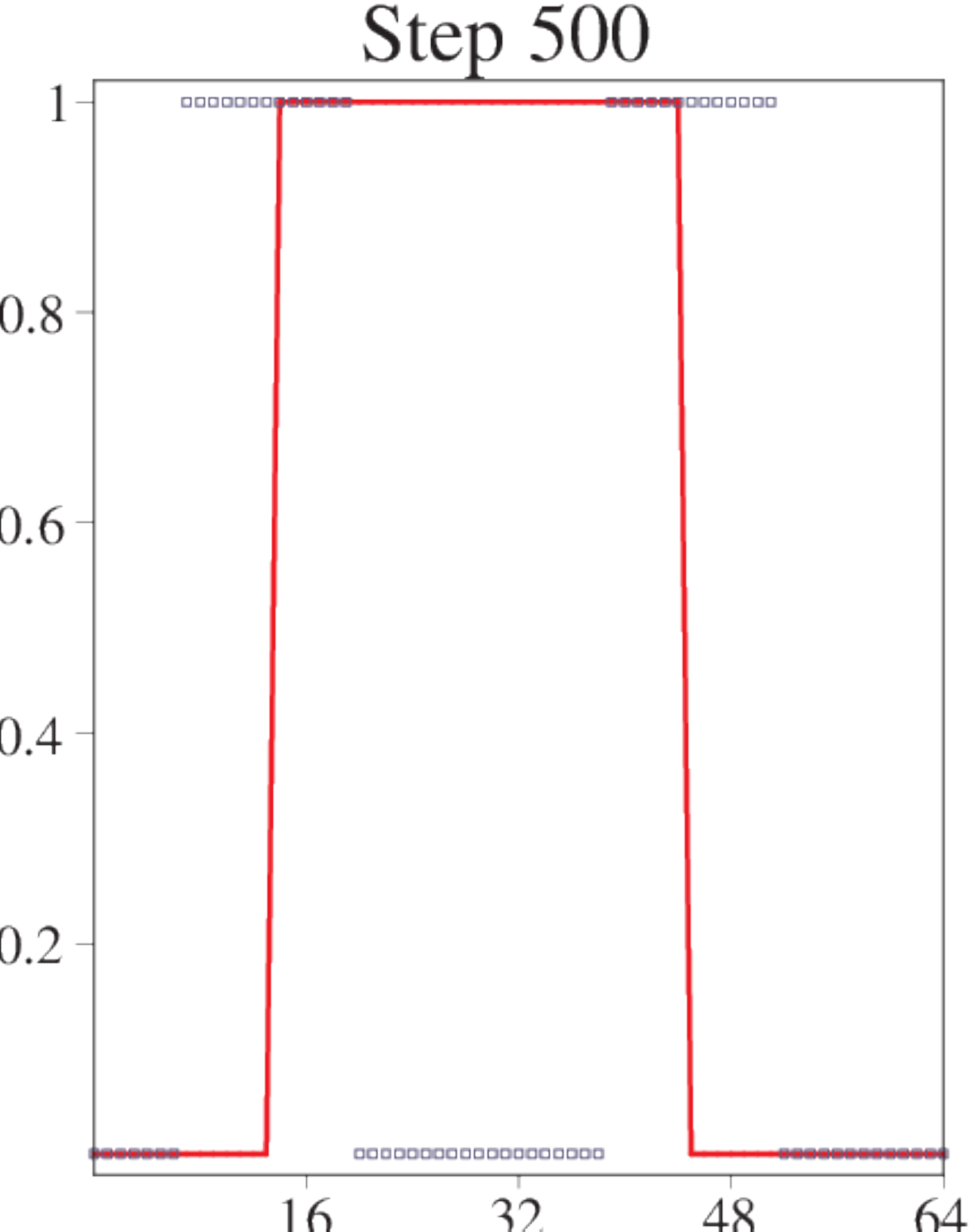} \hfil
             \epsfxsize4.2cm \epsfysize3.6cm
             \epsfbox{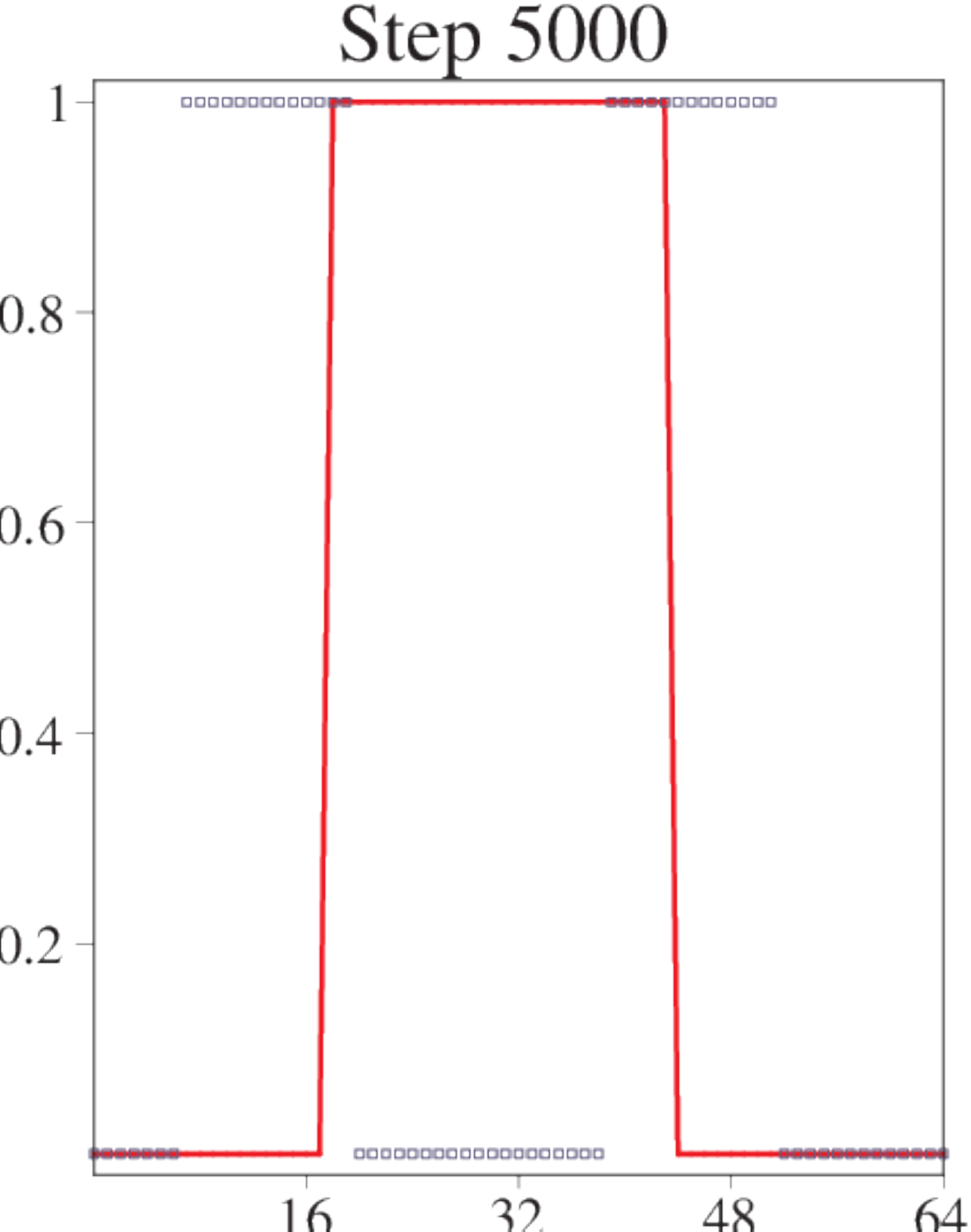} } \medskip
\centerline{ \epsfxsize4.2cm \epsfysize3.6cm
             \epsfbox{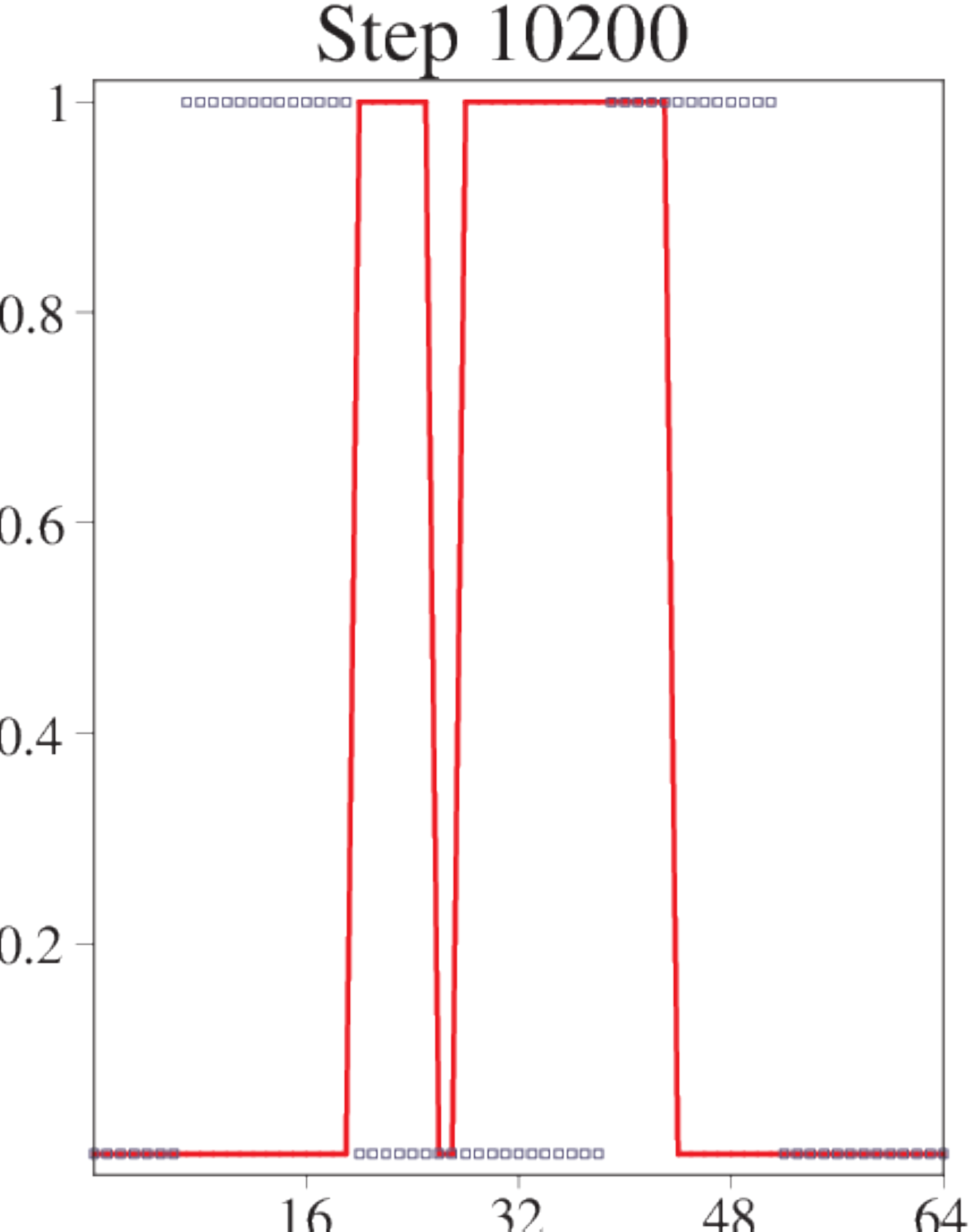} \hfil
             \epsfxsize4.2cm \epsfysize3.6cm
             \epsfbox{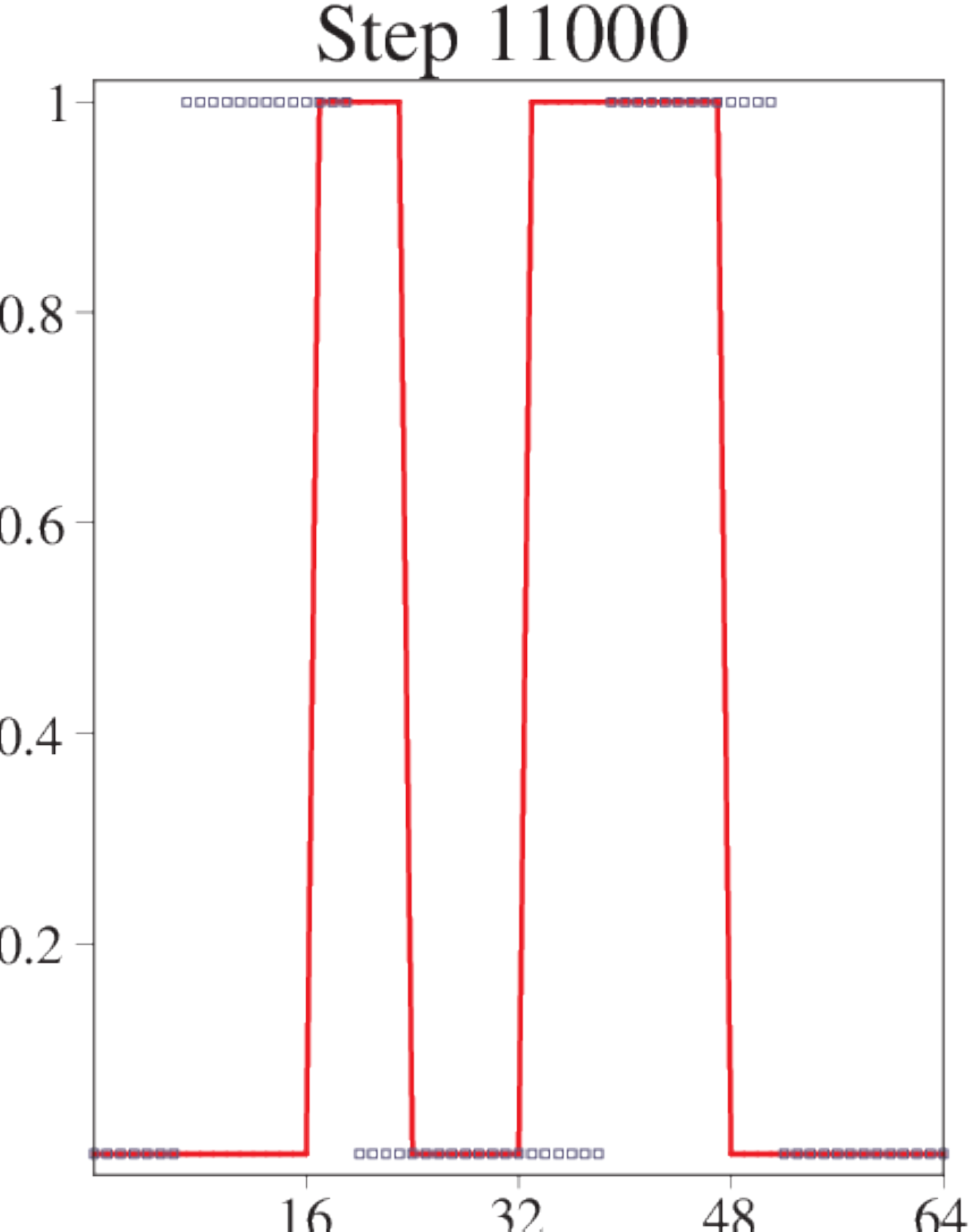} \hfil
             \epsfxsize4.2cm \epsfysize3.6cm
             \epsfbox{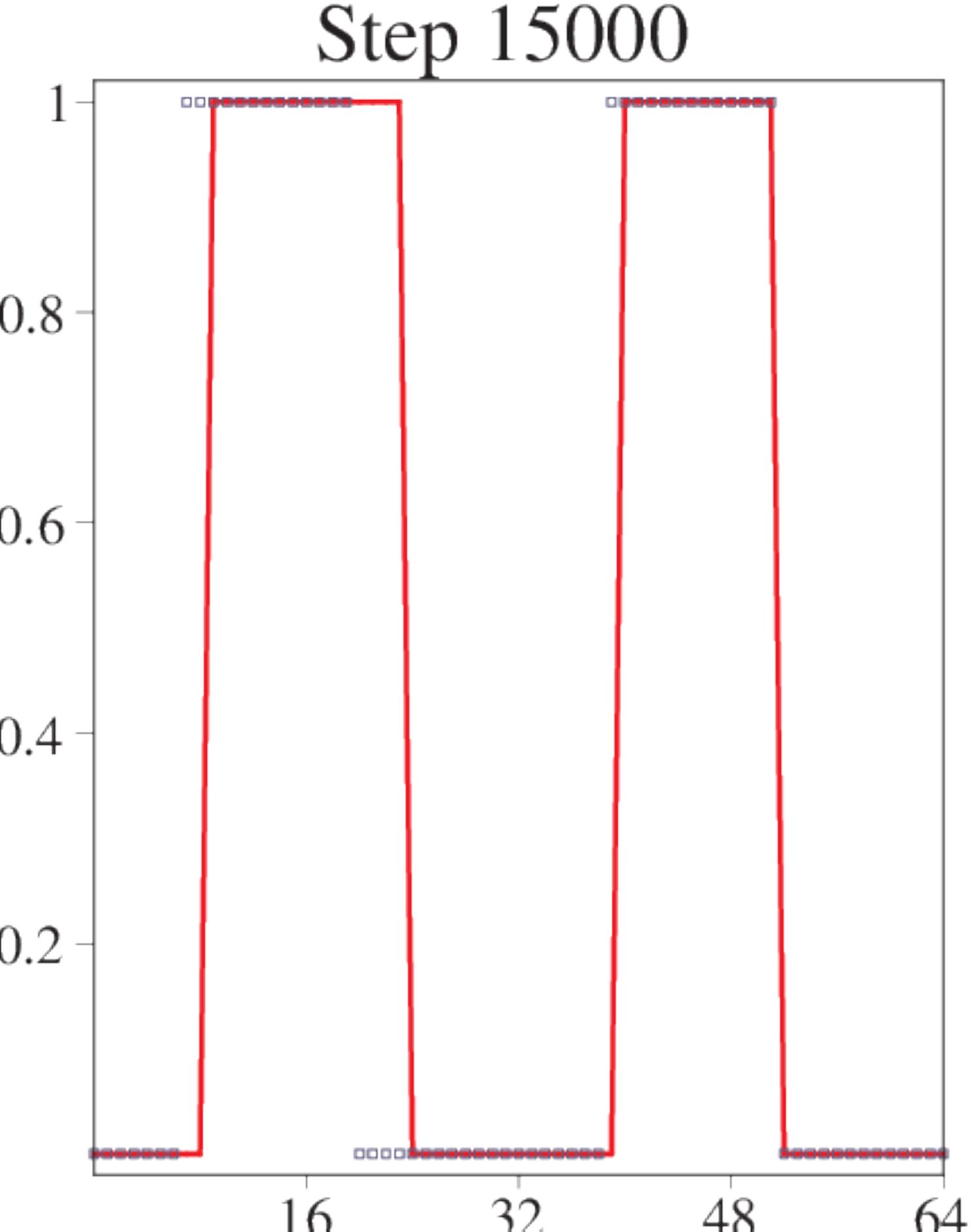} }
\caption{\footnotesize First numerical experiment:
On the top left, the initial guess for the level set function (solid red line).
The other pictures show the evolution of the level set method (solid red
line) after 500, 5000, 10200, 11000 and 15000 iterative steps. The the dotted
blue line represents the exact solution.
\label{fig:exp1-evol-ls}}
\end{figure}

\begin{figure}[th]
\centerline{ \epsfxsize4.2cm \epsfysize3.6cm
             \epsfbox{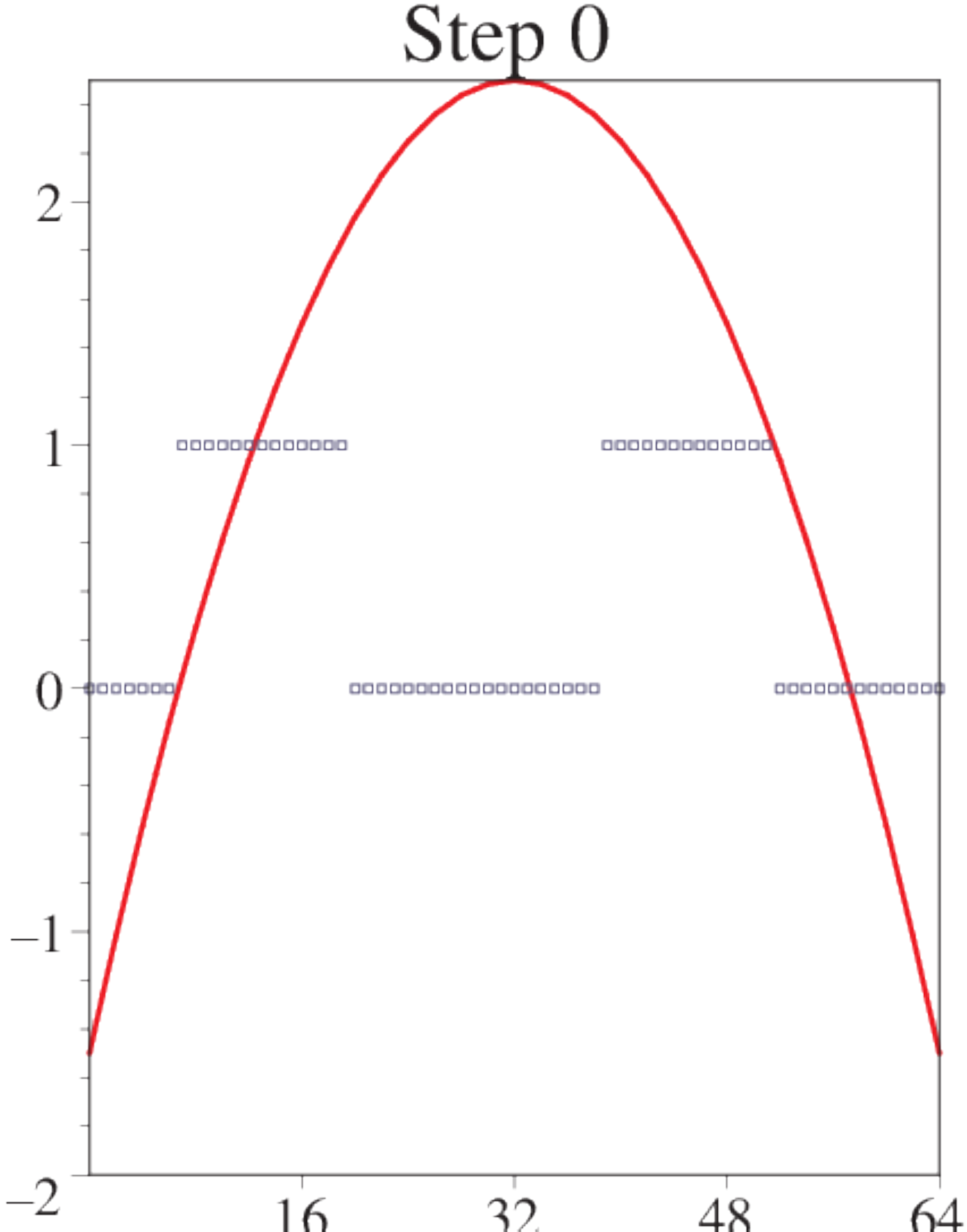} \hfil
             \epsfxsize4.2cm \epsfysize3.6cm
             \epsfbox{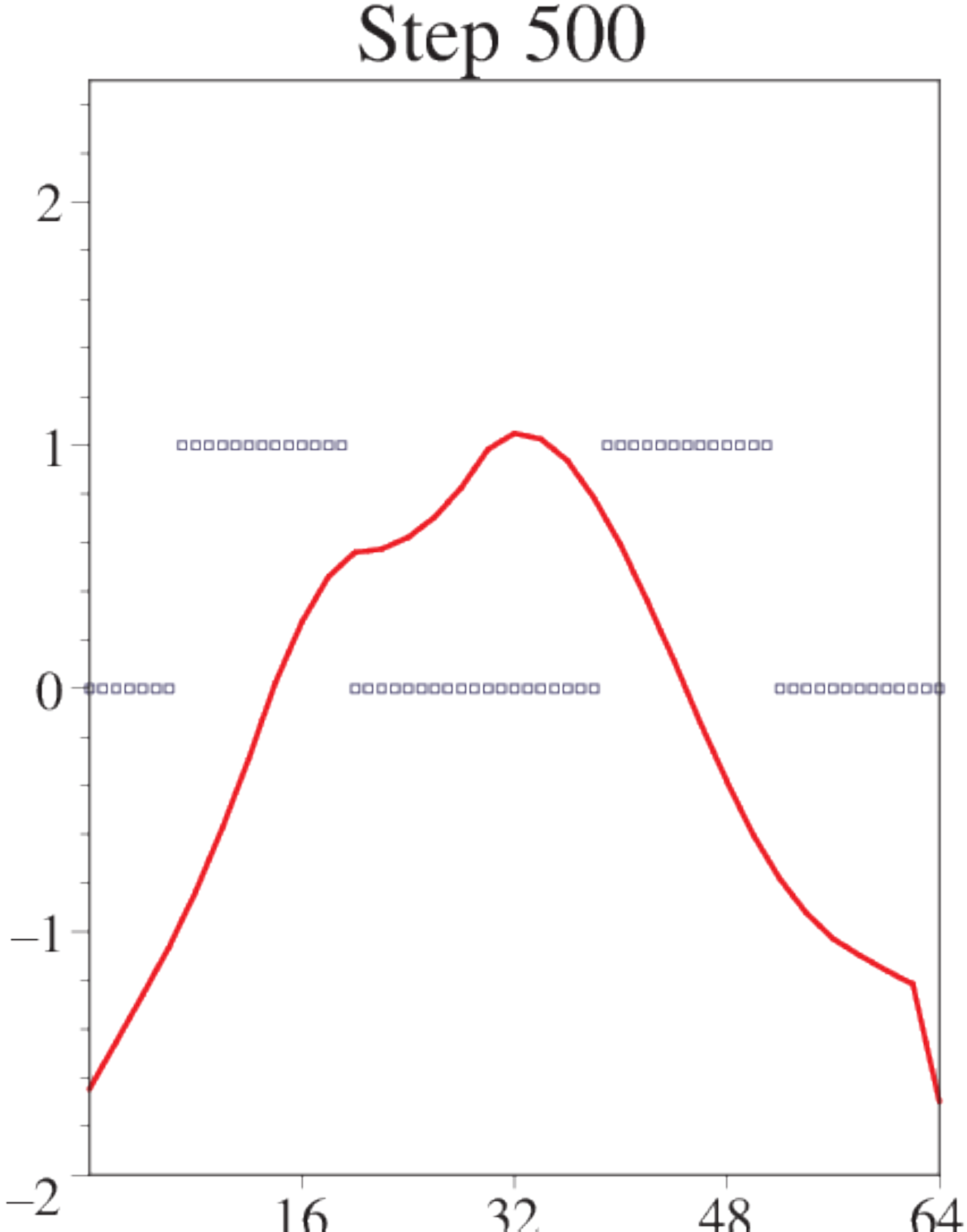} \hfil
             \epsfxsize4.2cm \epsfysize3.6cm
             \epsfbox{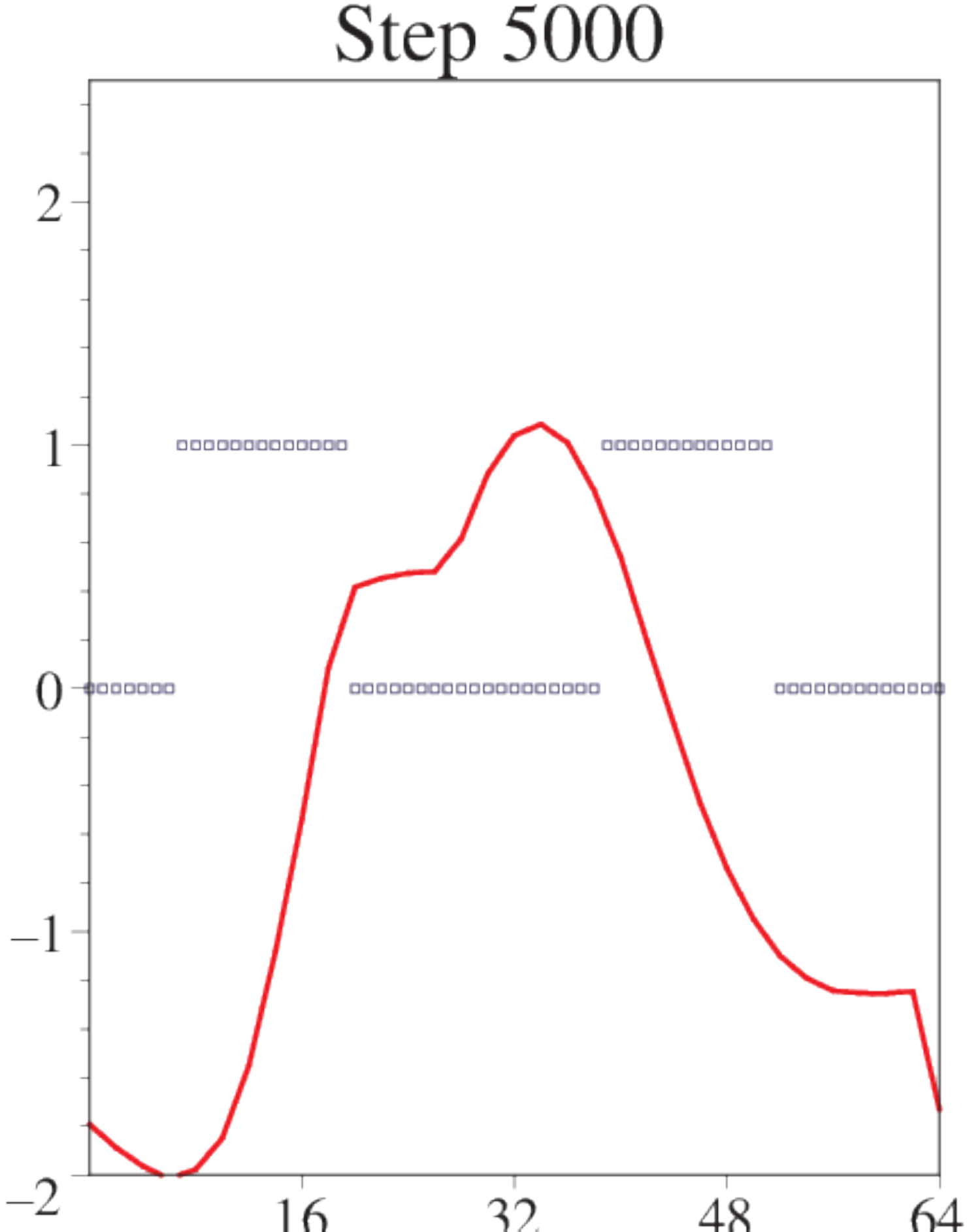} } \medskip
\centerline{ \epsfxsize4.2cm \epsfysize3.6cm
             \epsfbox{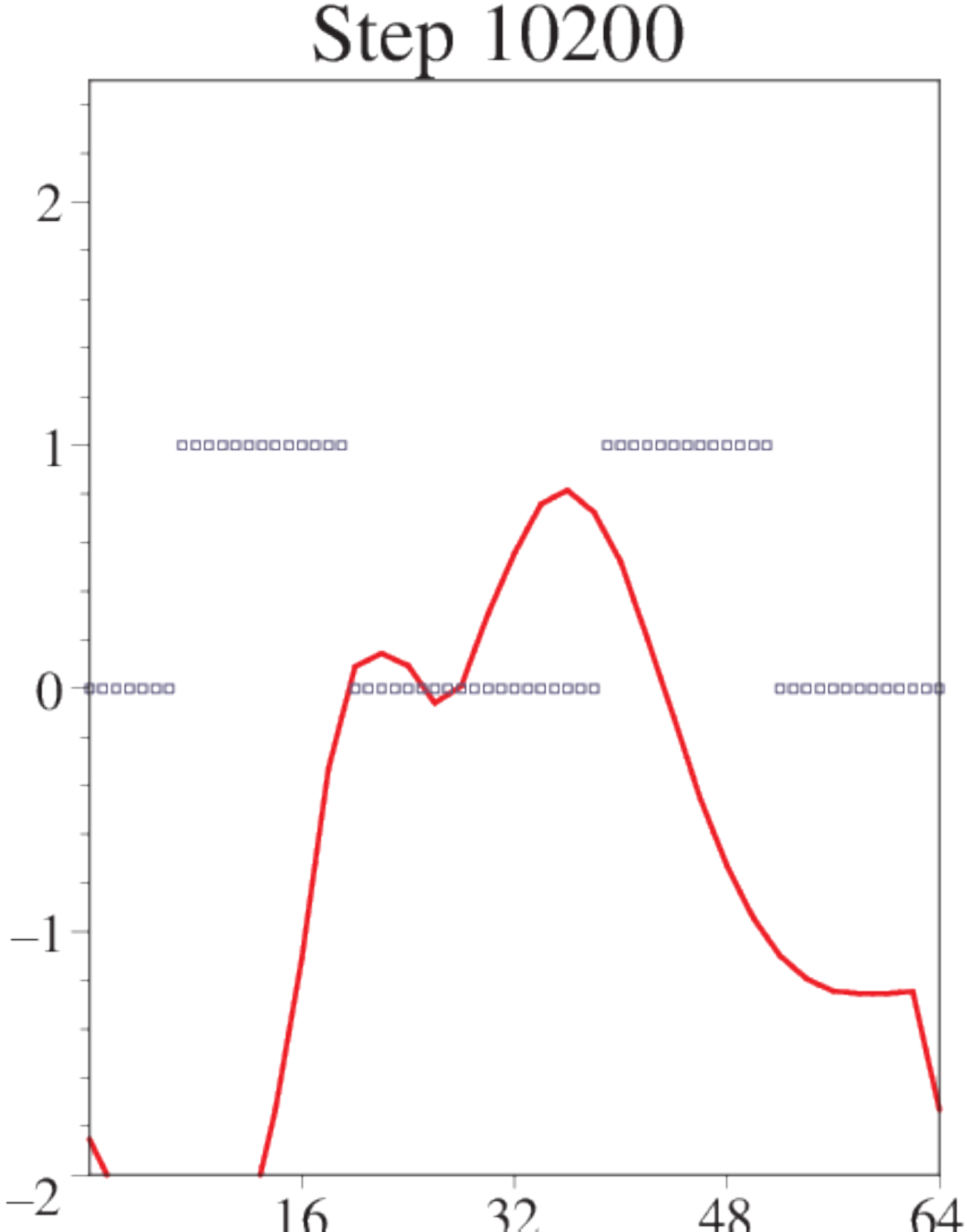} \hfil
             \epsfxsize4.2cm \epsfysize3.6cm
             \epsfbox{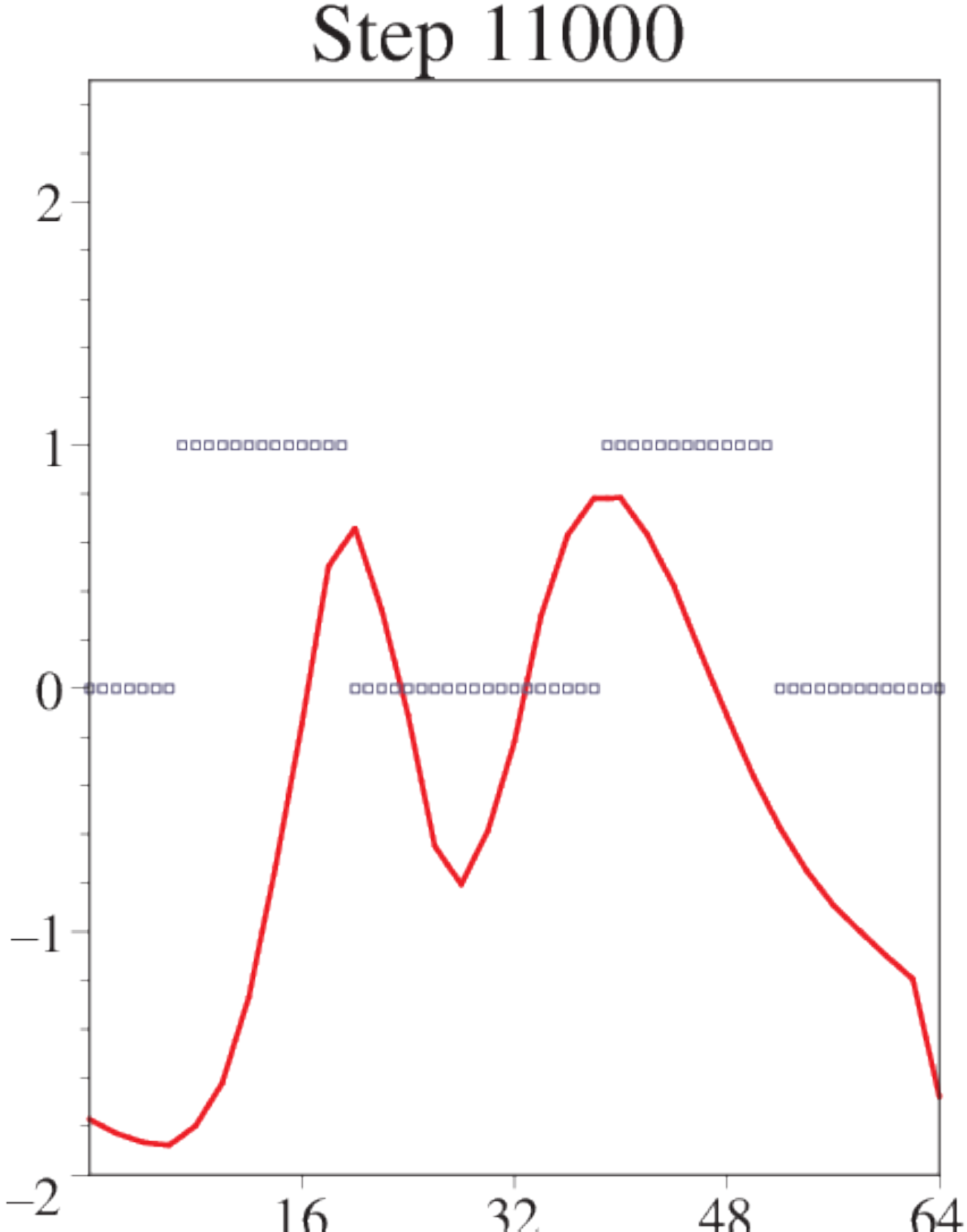} \hfil
             \epsfxsize4.2cm \epsfysize3.6cm
             \epsfbox{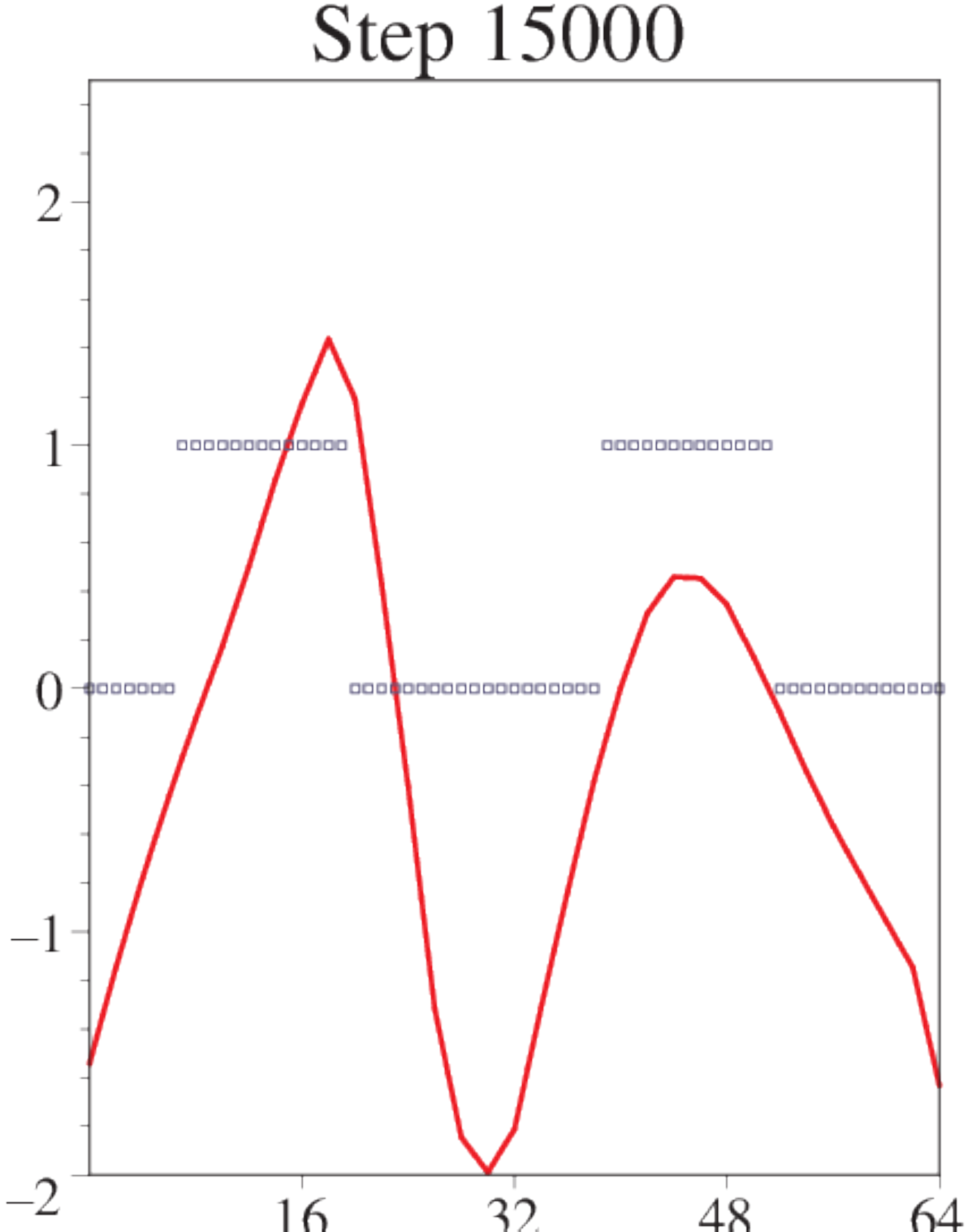} }
\caption{\footnotesize First numerical experiment:
On the top left, the initial guess for the level set function (solid red line).
The other pictures show the evolution of the level set function (solid red
line) after 500, 5000, 10200, 11000 and 15000 iterative steps. In all pictures
the exact solution $\overline{\vphi}$ of the Cauchy problem is represented by
the dotted blue line. \label{fig:exp1-evol-fk}}
\end{figure}

In Figure~\ref{fig:exp1-evol-ls} we show the level set iterations for the
Cauchy problem (\ref{eq:cp-exp1}) and in Figure~\ref{fig:exp1-evol-fk} the
corresponding evolution of level set function. Notice the large number of
steps required to obtain a precise approximation. As already observed when
applying level set methods to other models \cite{FSL05,San96}, the splitting
of the level sets happens only after a large number of iterations (over
10000 in this experiment). Nevertheless, the {\em a priori} required
precision could always be reached (in this experiment, 
$\| \vphi_k - \overline{\vphi} \|_{L^2(\Gamma_2)} < 10^{-2}$ was reached
after 15000 steps).
It is worth noticing that, in our experiments, the total number of steps
needed to reach a pre-specified precision does not depend from the initial
guess for the level set method.

\subsection{Degree of ill-posedness affecting convergence rates}
\label{ssec:5-2}
\begin{figure}[ht]
\centerline{ \epsfxsize6.3cm \epsfbox{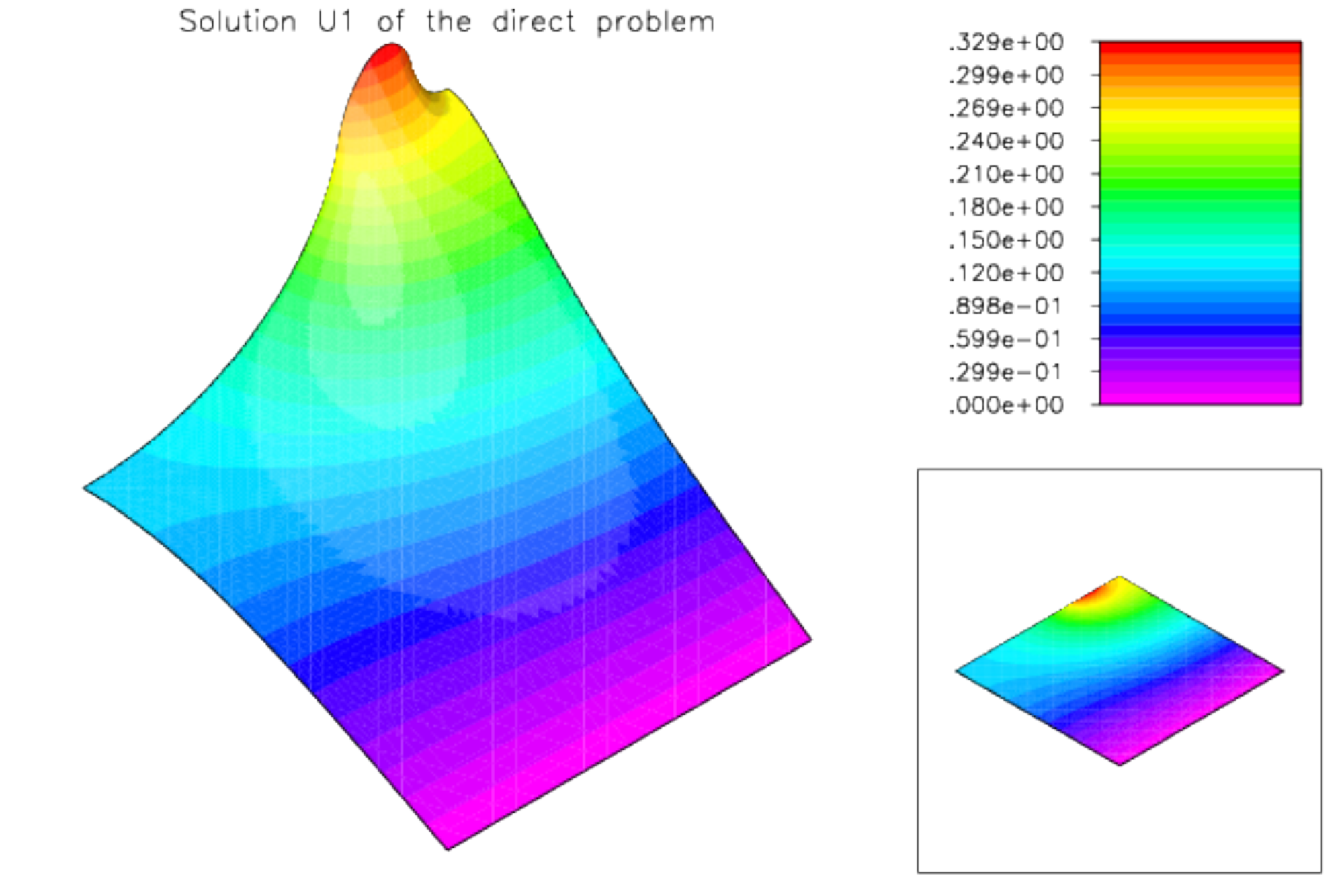} \hfil
             \epsfxsize6.3cm \epsfbox{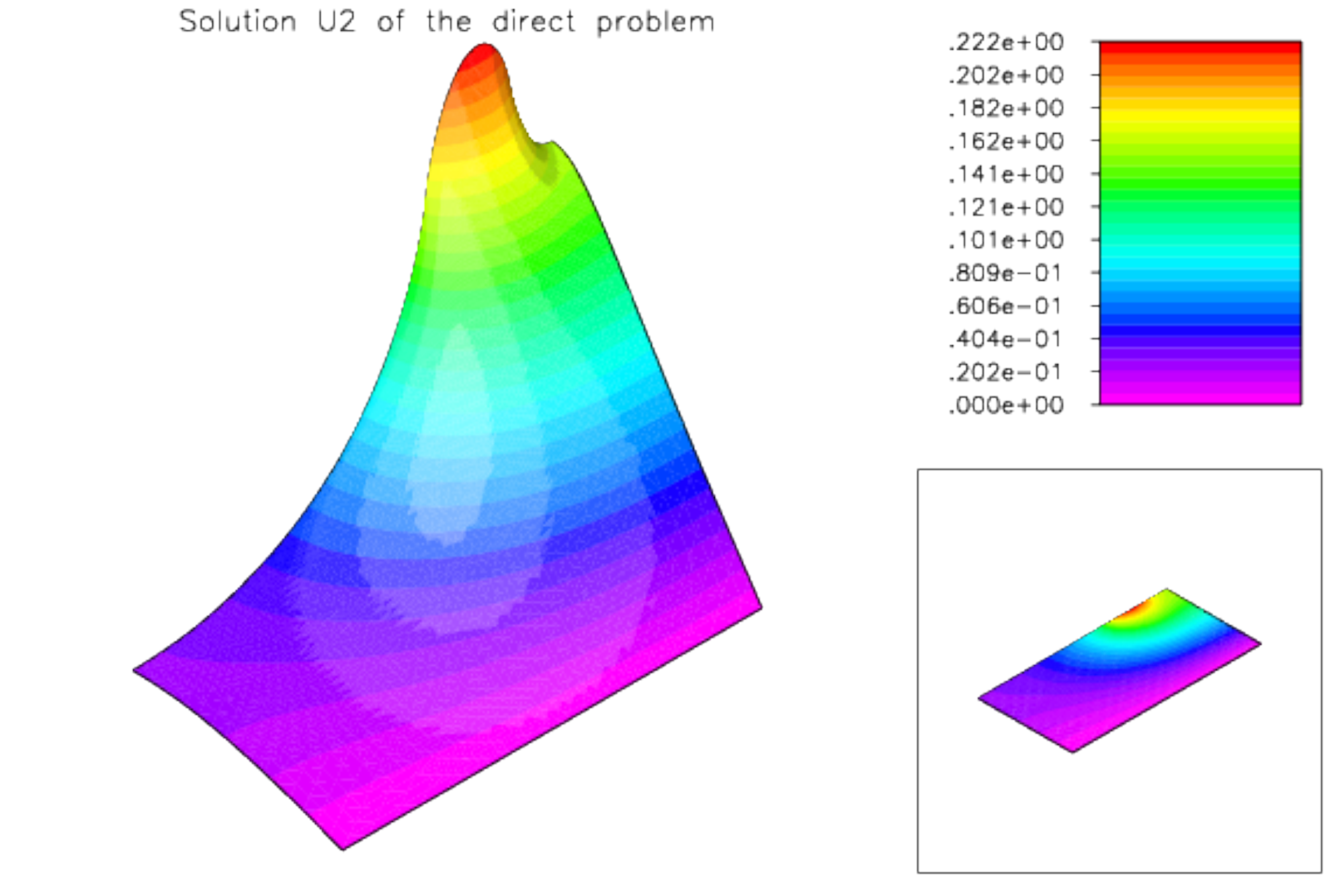}}
\caption{\footnotesize Framework for the second numerical experiment:
The solutions $u_i$ of the elliptic BVPs corresponding to (CP$)_i$, $i=1,2$,
are depicted on the left and on the right hand side respectively. The Cauchy
data $g_{2,i}$ correspond to $(u_i)_\nu$ at $\Gamma_1$ (lower edge). Notice
also that the Dirichlet boundary condition $u_1 = u_2 = 0$ at $\Gamma_1$
holds.
Both Cauchy problems share the same solution $(u_1)_\nu |_{\Gamma_{2,1}} =
\overline{\vphi} = (u_2)_\nu |_{\Gamma_{2,2}}$ (top edge).
\label{fig:exp2-frwork}}
\end{figure}

In this experiment we analyze the effect of the ill-posedness degree of
the Cauchy problem on the performance of the level set method introduced in
Section~\ref{sec:4}. For this purpose we introduce the domains
$\Omega_1 := (0,1) \times (0,1)$ and $\Omega_2 := (0,1) \times (0,0.5)$.
The boundary part $\Gamma_1 := (0,1) \times \{0\}$ is the same for both
Domains. Moreover we define $\Gamma_{2,1} := (0,1) \times \{1\}$, 
$\Gamma_{2,2} := (0,1) \times \{0.5\}$ and $\Gamma_{3,i} := \partial\Omega_i
/ \{ \Gamma_1 \cup \Gamma_{2,i}\}$, $i=1,2$.

\begin{figure}[th]
\centerline{ \epsfxsize4.2cm \epsfysize3.6cm
             \epsfbox{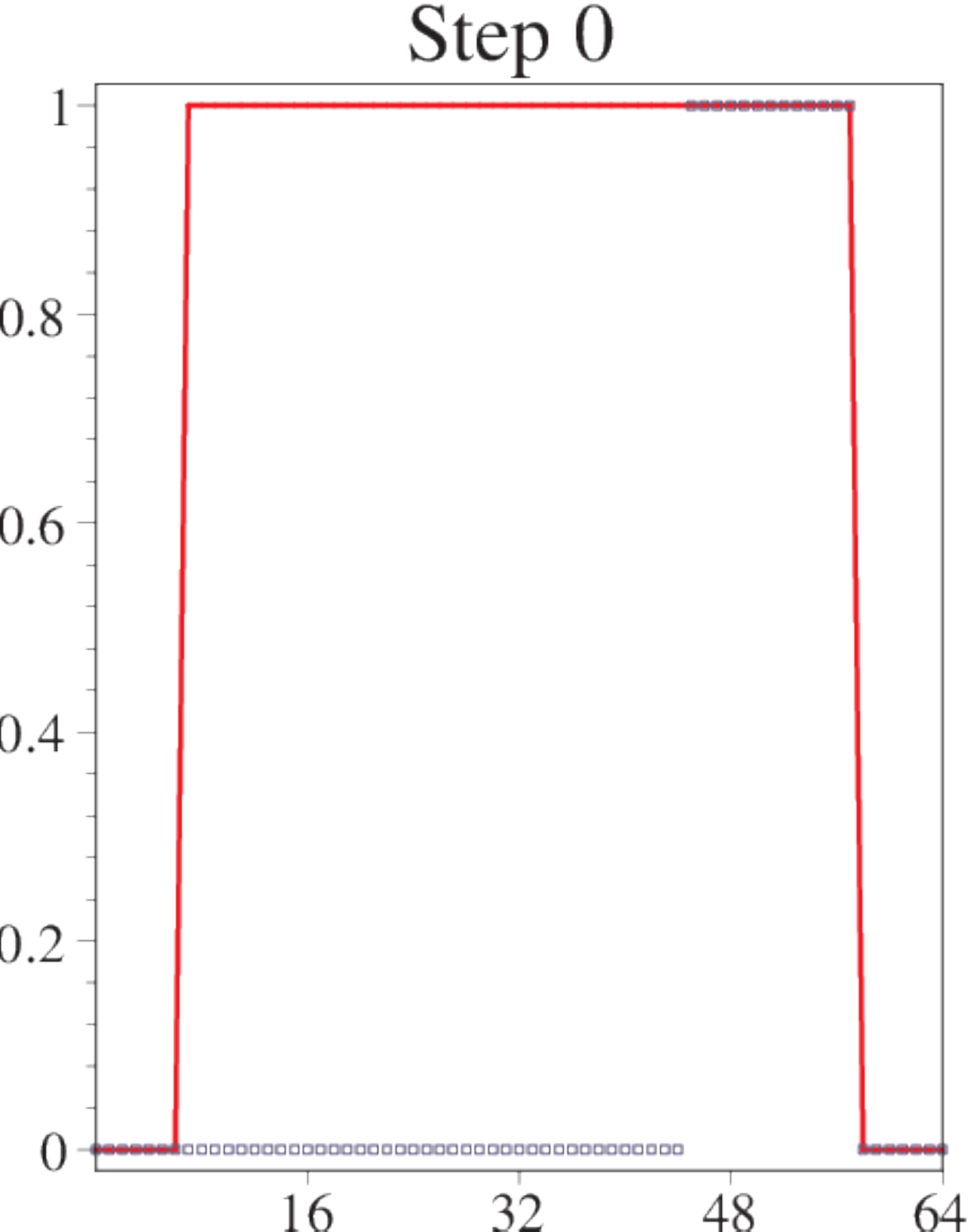} \hfil
             \epsfxsize4.2cm \epsfysize3.6cm
             \epsfbox{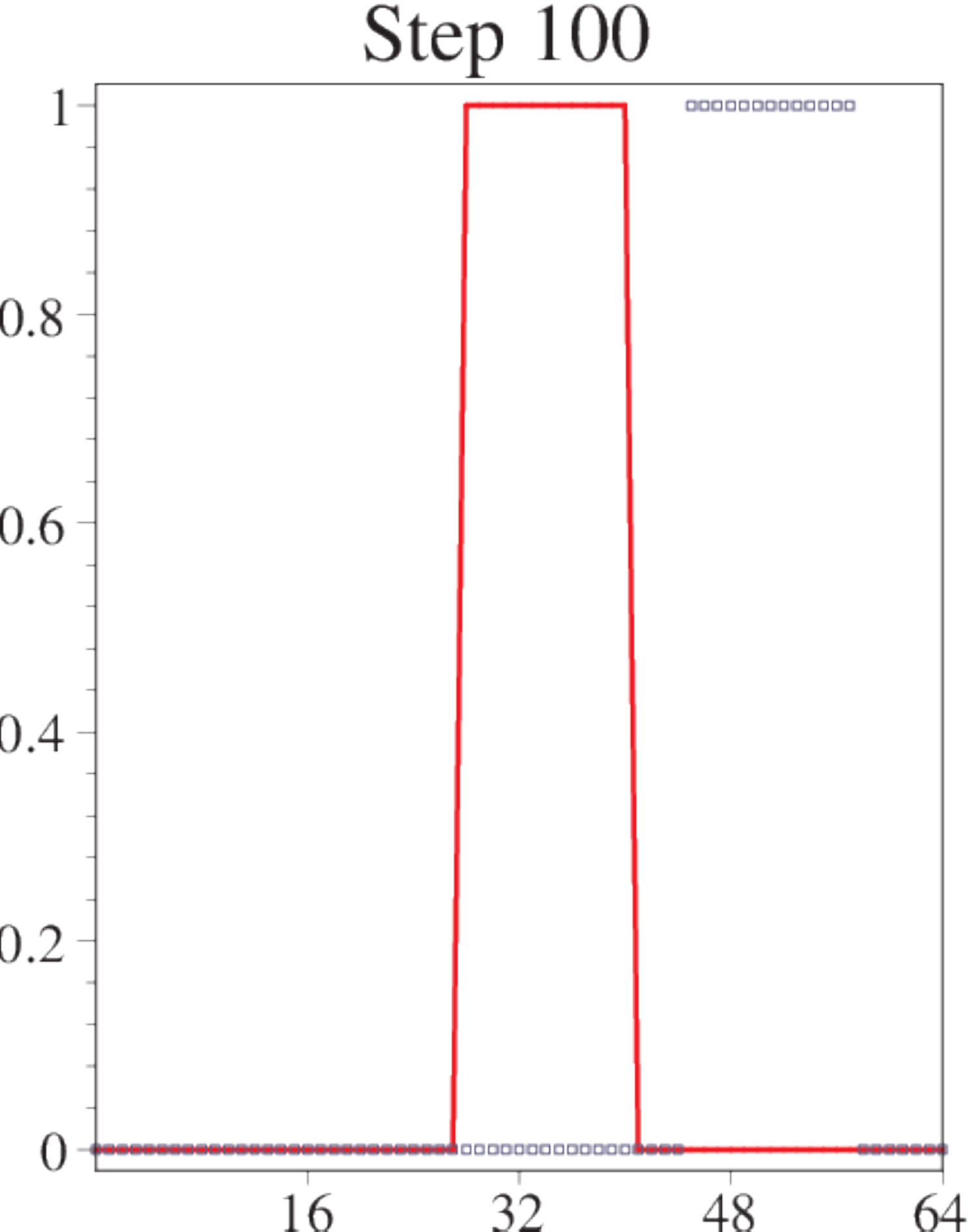} \hfil
             \epsfxsize4.2cm \epsfysize3.6cm
             \epsfbox{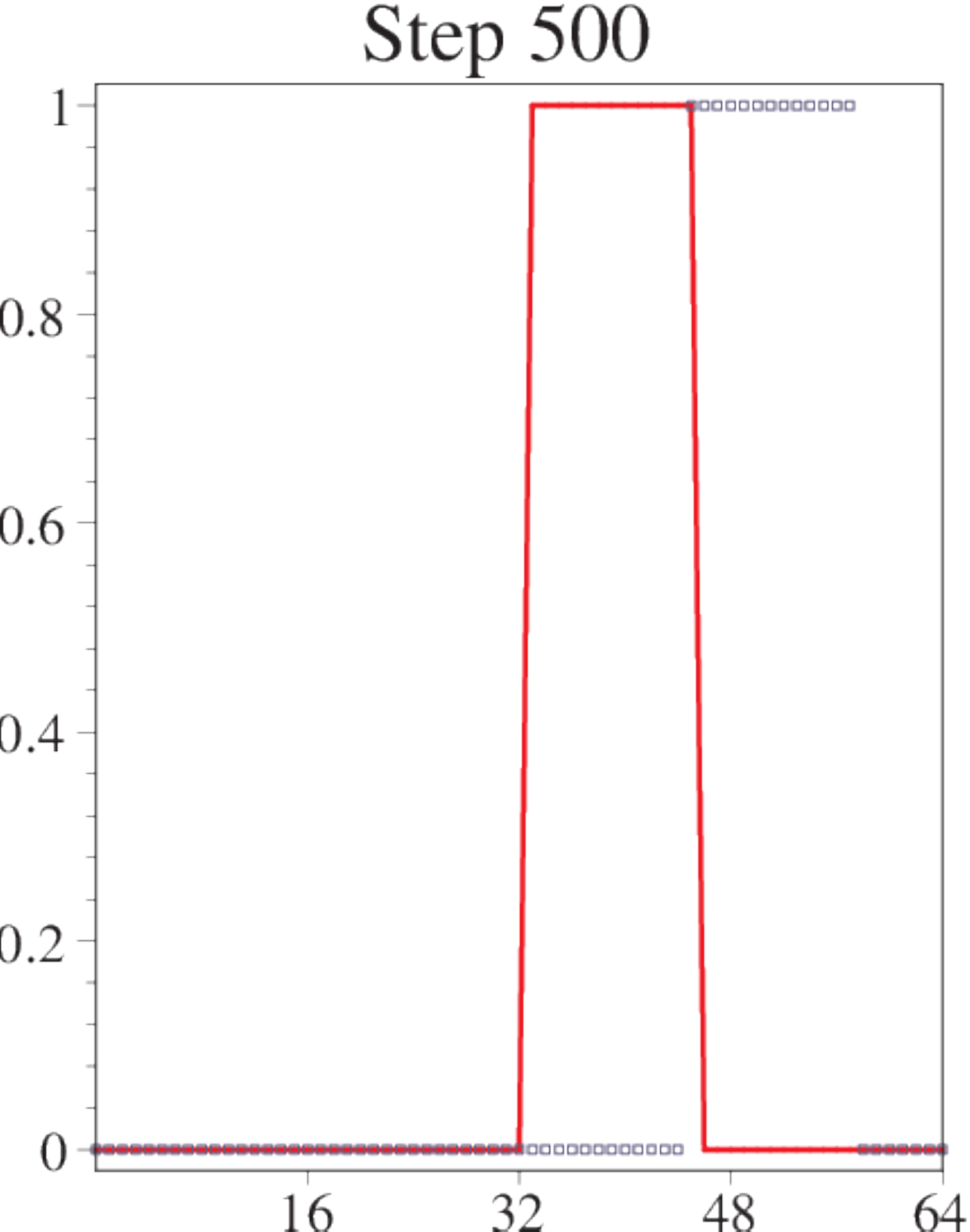} } \medskip
\centerline{ \epsfxsize4.2cm \epsfysize3.6cm
             \epsfbox{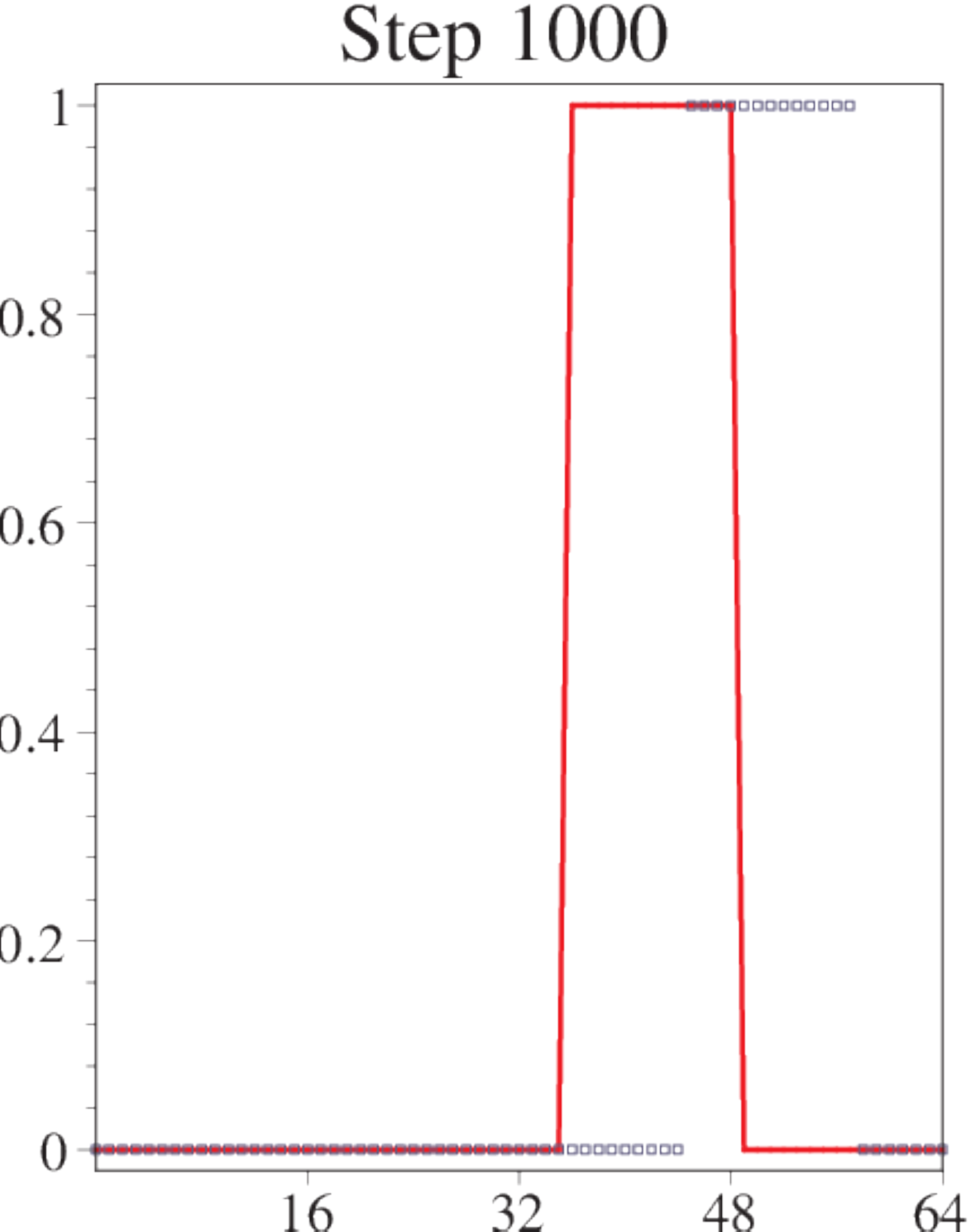} \hfil
             \epsfxsize4.2cm \epsfysize3.6cm
             \epsfbox{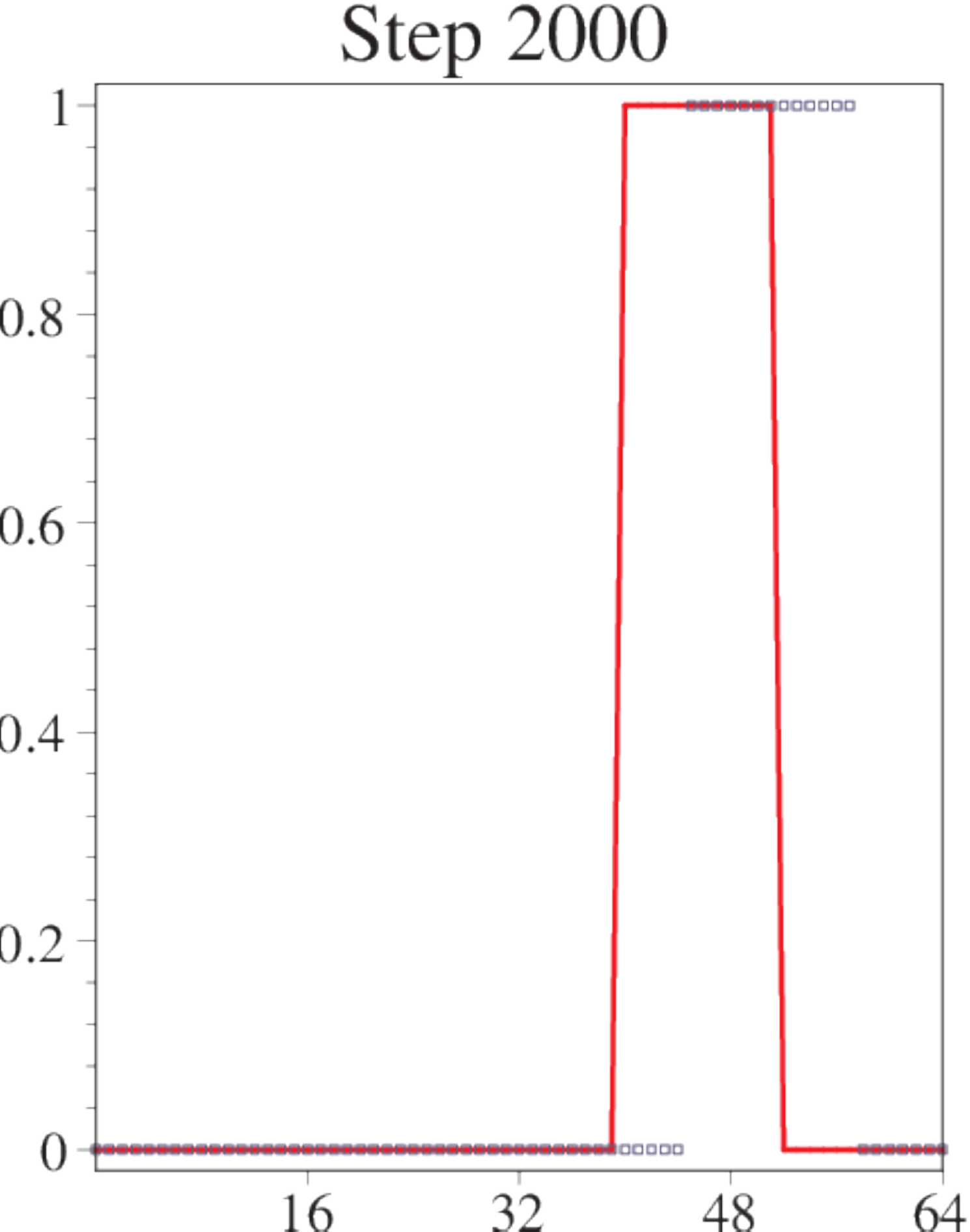} \hfil
             \epsfxsize4.2cm \epsfysize3.6cm
             \epsfbox{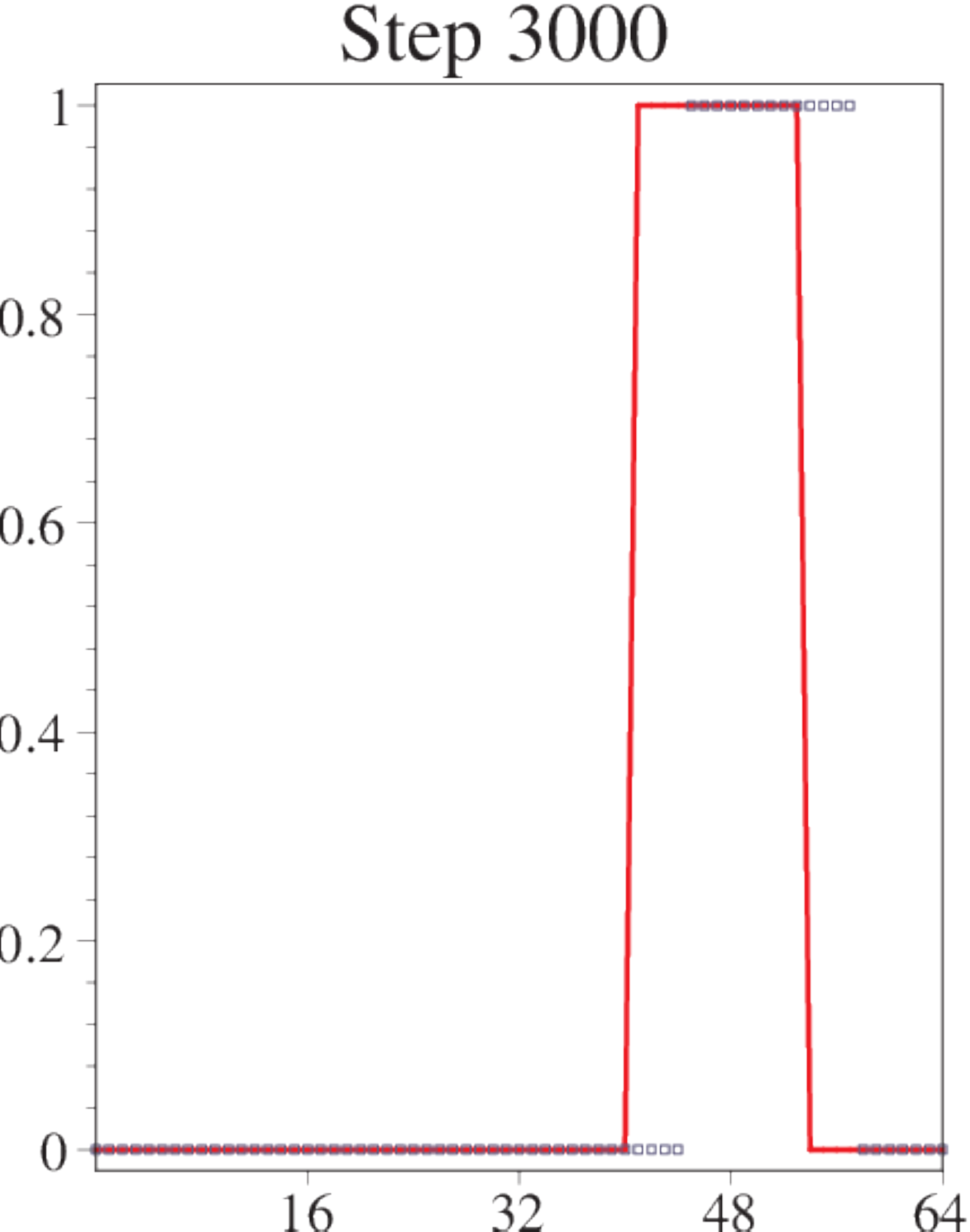} }
\caption{\footnotesize Second numerical experiment:
On the top left, the exact solution $\overline{\vphi}$ of the Cauchy problems
(dotted blue line) and the initial guess for the level set method (solid red
line).
The other pictures show the evolution of the level set method for problem
(CP$)_1$ after 100, 500, 1000, 2000 and 3000 iterative steps. The solid
(red) line represent the iteration and the dotted (blue) line the exact
solution.\label{fig:exp2-evol-ls1}}
\end{figure}

We compare the performance of our method for the Cauchy problems (CP$)_i$,
$i=1, 2$, defined by
$$  \Delta u_i  =  0      \, , \ \mbox{in } \, \Omega_i \quad
    u_i  =  0             \, , \ \mbox{at } \, \Gamma_1 \quad
    (u_i)_\nu  =  g_{2,i} \, , \ \mbox{at } \, \Gamma_1 \quad
    (u_i)_\nu  =  0       \, , \ \mbox{at } \, \Gamma_3 \, . $$
The Cauchy data $(0,g_{2,i})$ is chosen such that problems (CP$)_i$ have the
same solution, i.e. $(u_1)_\nu |_{\Gamma_{2,1}} = (u_2)_\nu |_{\Gamma_{2,2}}$
(see Figures~\ref{fig:exp2-frwork}, \ref{fig:exp2-evol-ls1}).

\begin{figure}[th]
\centerline{ \epsfxsize4.2cm \epsfysize3.6cm
             \epsfbox{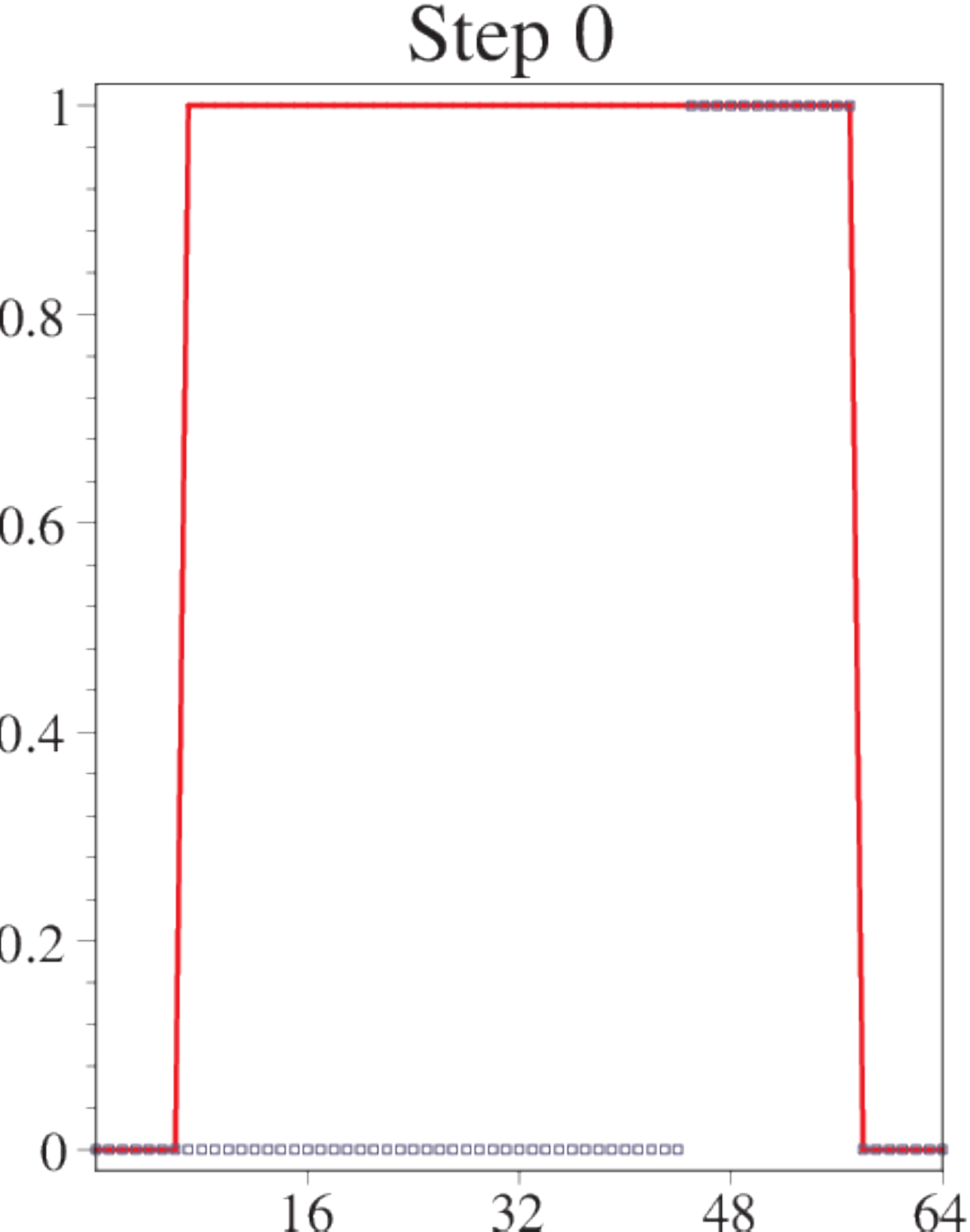} \hfil
             \epsfxsize4.2cm \epsfysize3.6cm
             \epsfbox{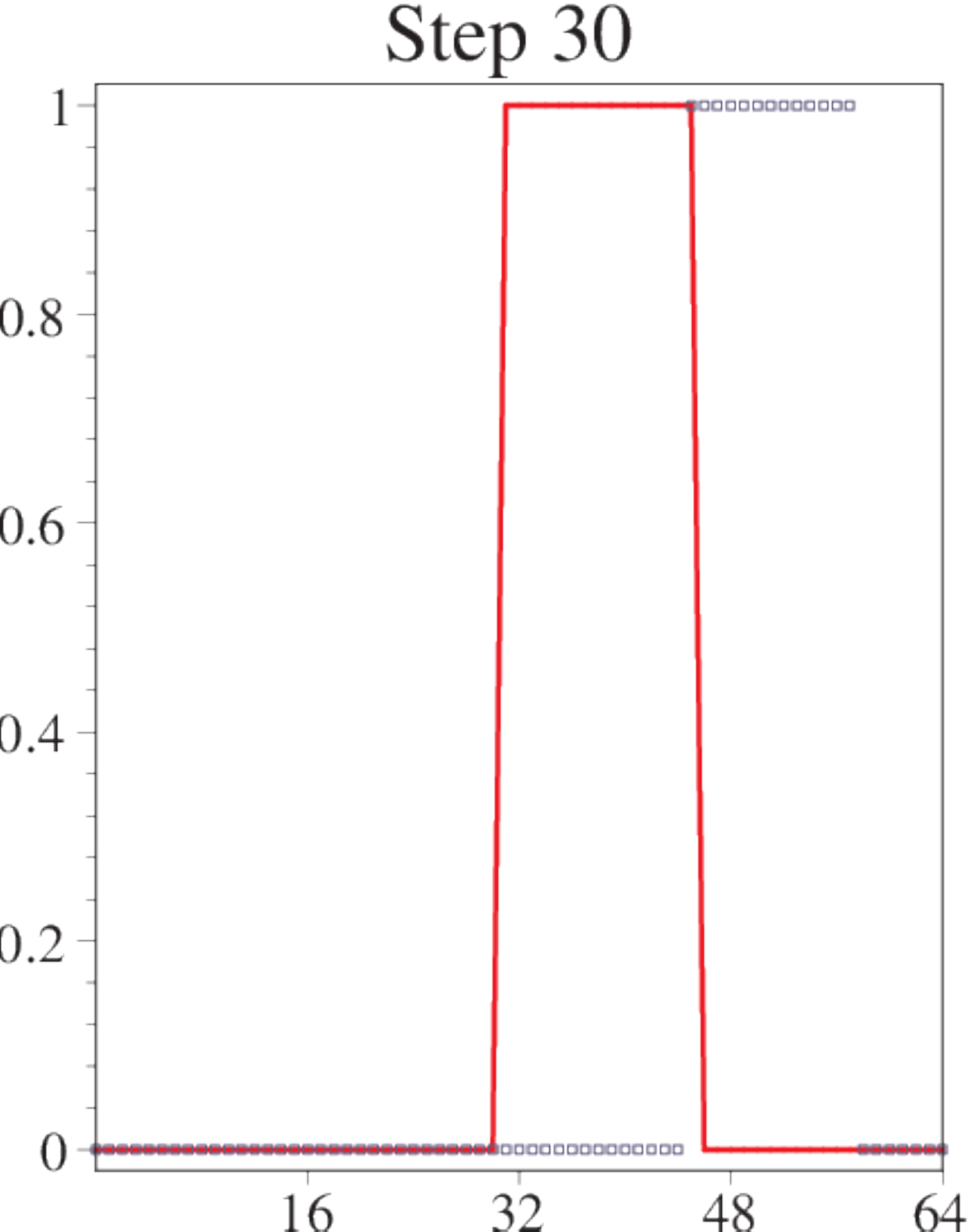} \hfil
             \epsfxsize4.2cm \epsfysize3.6cm
             \epsfbox{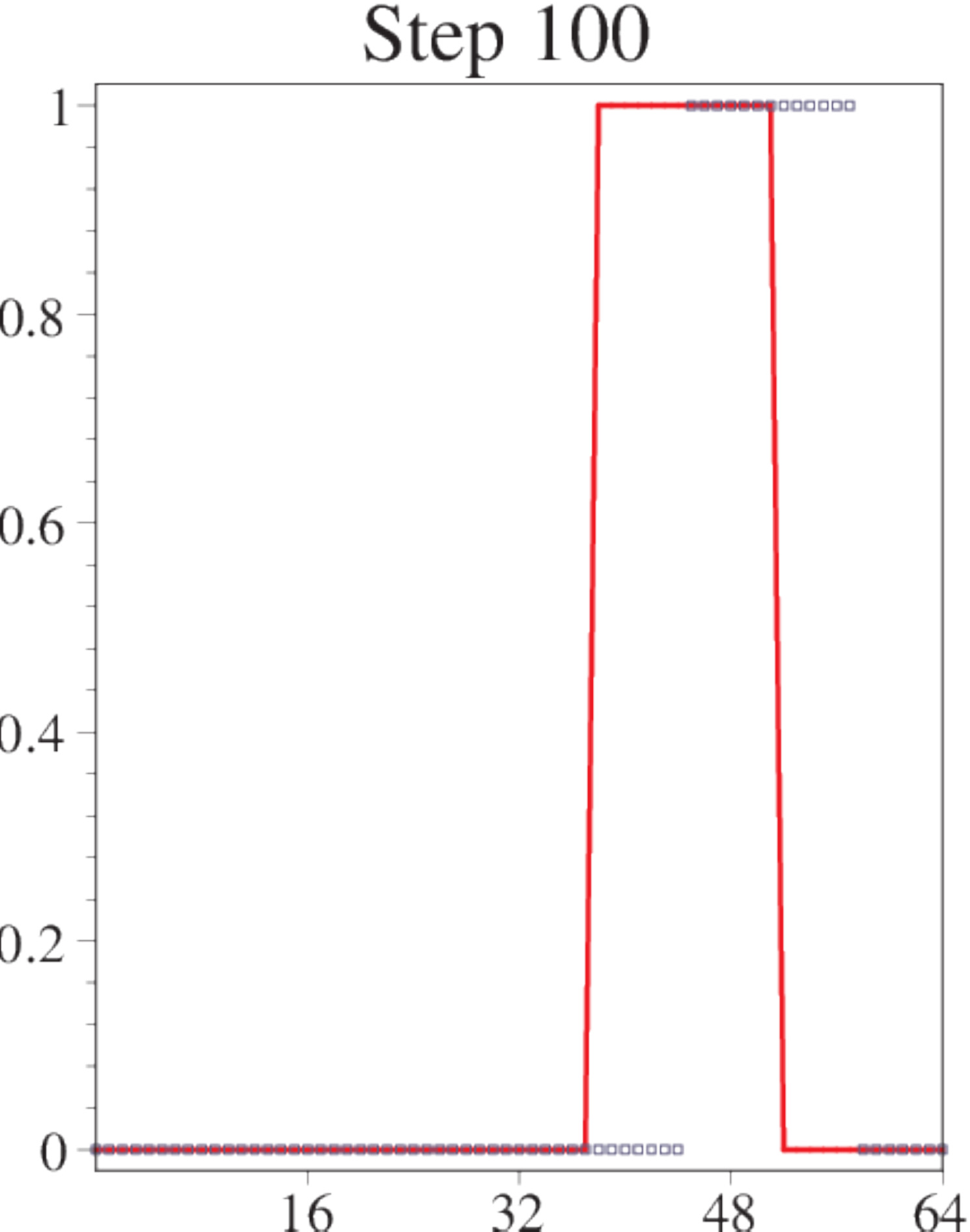} } \medskip
\centerline{ \epsfxsize4.2cm \epsfysize3.6cm
             \epsfbox{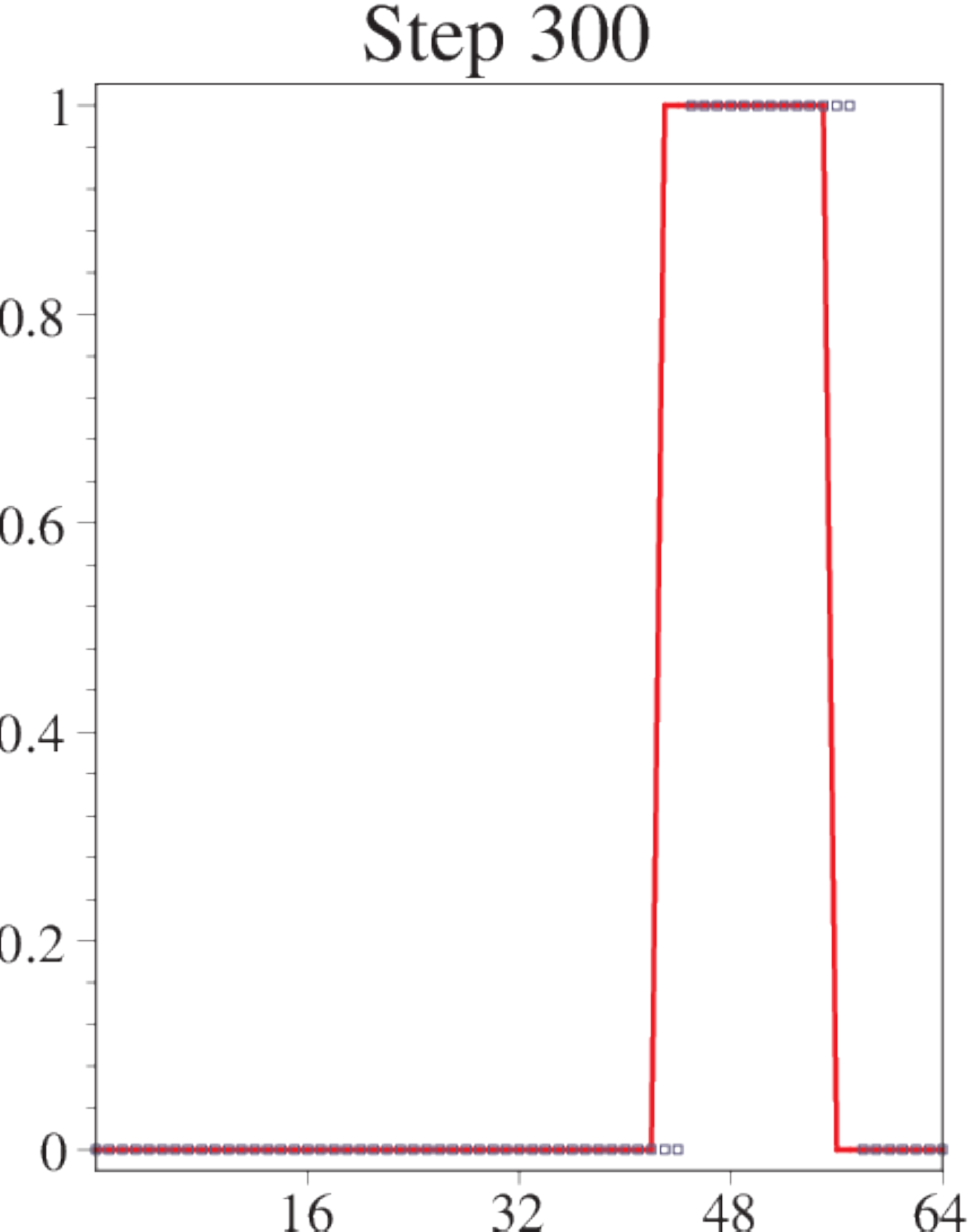} \hfil
             \epsfxsize4.2cm \epsfysize3.6cm
             \epsfbox{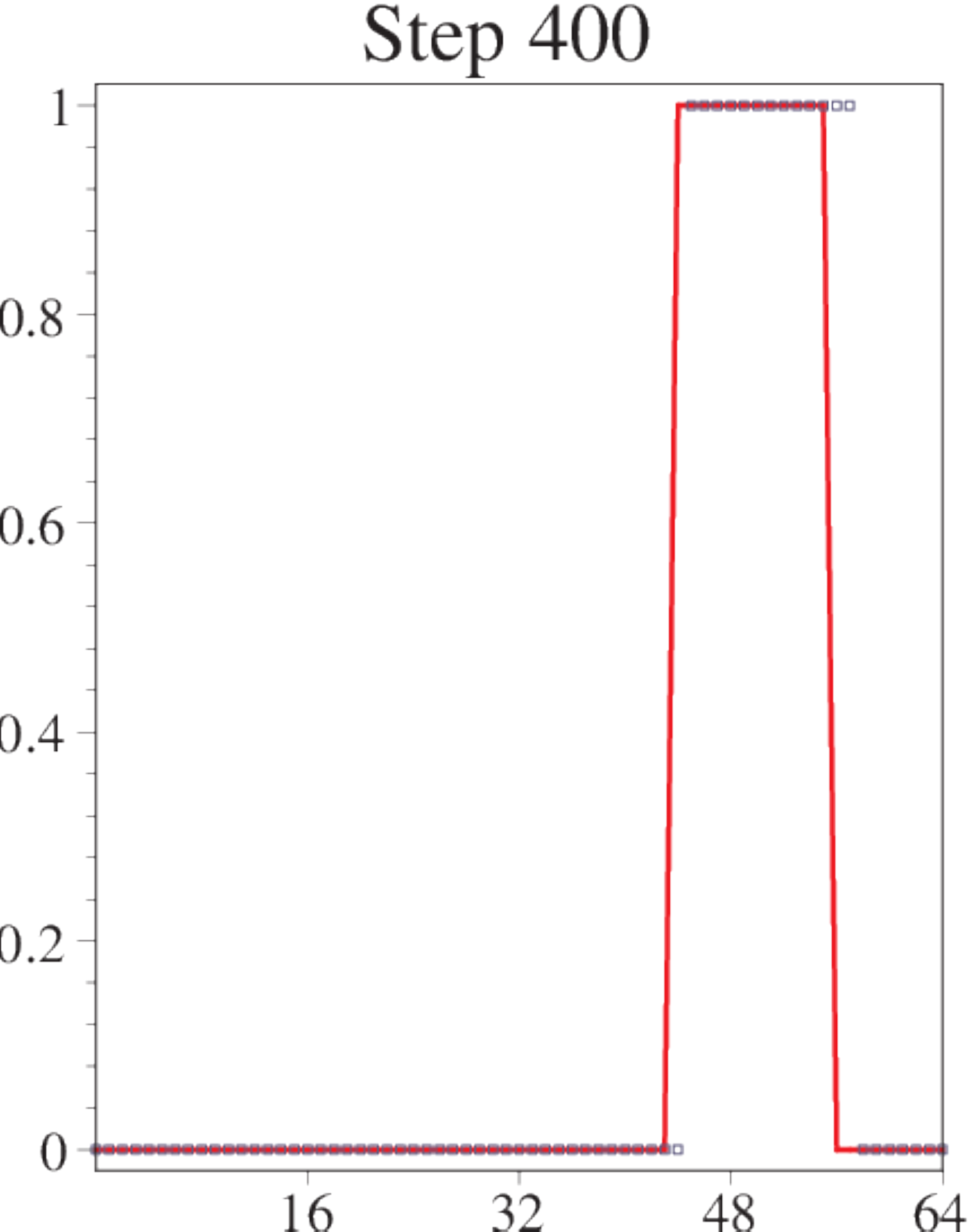} \hfil
             \epsfxsize4.2cm \epsfysize3.6cm
             \epsfbox{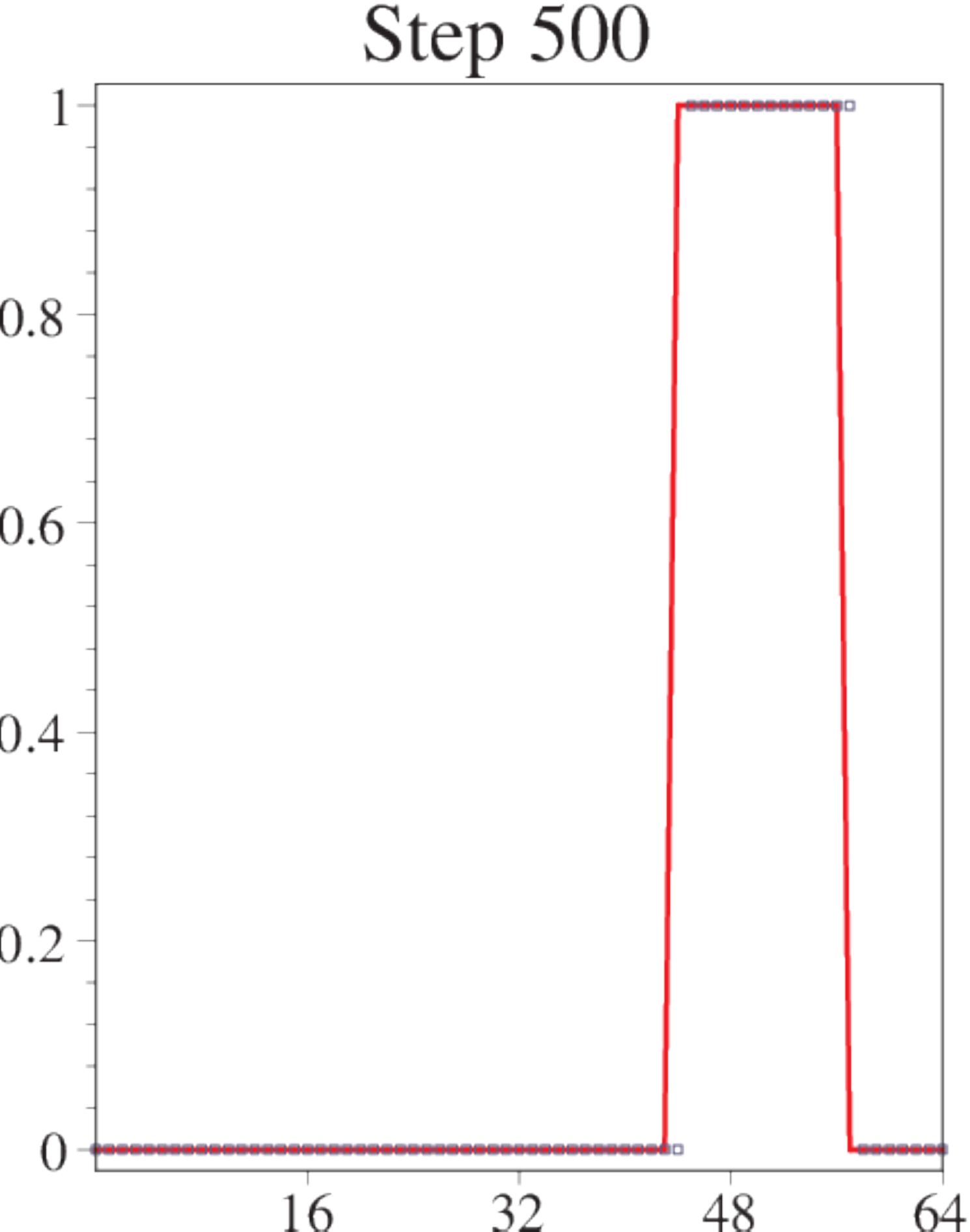} }
\caption{\footnotesize Second numerical experiment:
Evolution of the level set method for problem (CP$)_2$. Plots of the initial
guess and after 30, 100, 300, 400 and 500 iterative steps. The solid
(red) line represent the iteration and the dotted (blue) line the exact
solution.\label{fig:exp2-evol-ls2}}
\end{figure}

It is known from the literature that the ill-posedness degree of (CP)
increases with the {\em distance} between the boundary parts $\Gamma_1$
and $\Gamma_2$ \cite{Le00}. Since we are considering simple domains
$\Omega_i$, it is possible to compute the eigenvalues $\{ \lbd_{i,j} \}_j$
of the operators $L_i$ in (\ref{def:oper-L}) corresponding to (CP$)_i$
$$ \lbd_{i,j} = \sinh(j/i)^{-1} , \ j = 1, 2, \dots \, , \ i = 1, 2 \, . $$
This gives us a measure how ill-posed (CP$)_1$ is when compared with (CP$)_2$.
In the sequel we compare the performance of the level set method for both
problems (exact Cauchy data is used). In Figures~\ref{fig:exp2-evol-ls1} and
\ref{fig:exp2-evol-ls2} we see the level set iterations for problems (CP$)_1$
and (CP$)_2$ respectively.
Due to the difference between the ill-posedness degree of both problems, the
method converges faster for problem (CP$)_2$. Nevertheless, the same accuracy
can be reached for both problems if one iterates long enough. We conclude
that the degree of ill-posedness only affects the amount of computational
effort needed to obtain an approximate solution with {\em a priori} defined
precision, and not the quality of the final approximation.
Our experiments indicate that the number of iteration steps needed by the
level set method (to reach a desired accuracy) increases exponentially with
the distance between $\Gamma_1$ and $\Gamma_2$, the same way the degree of
ill-posedness also does.

\subsection{An experiment with noisy Cauchy data} \label{ssec:5-3}

\begin{figure}[th]
\centerline{ \epsfxsize4.2cm \epsfysize3.6cm
             \epsfbox{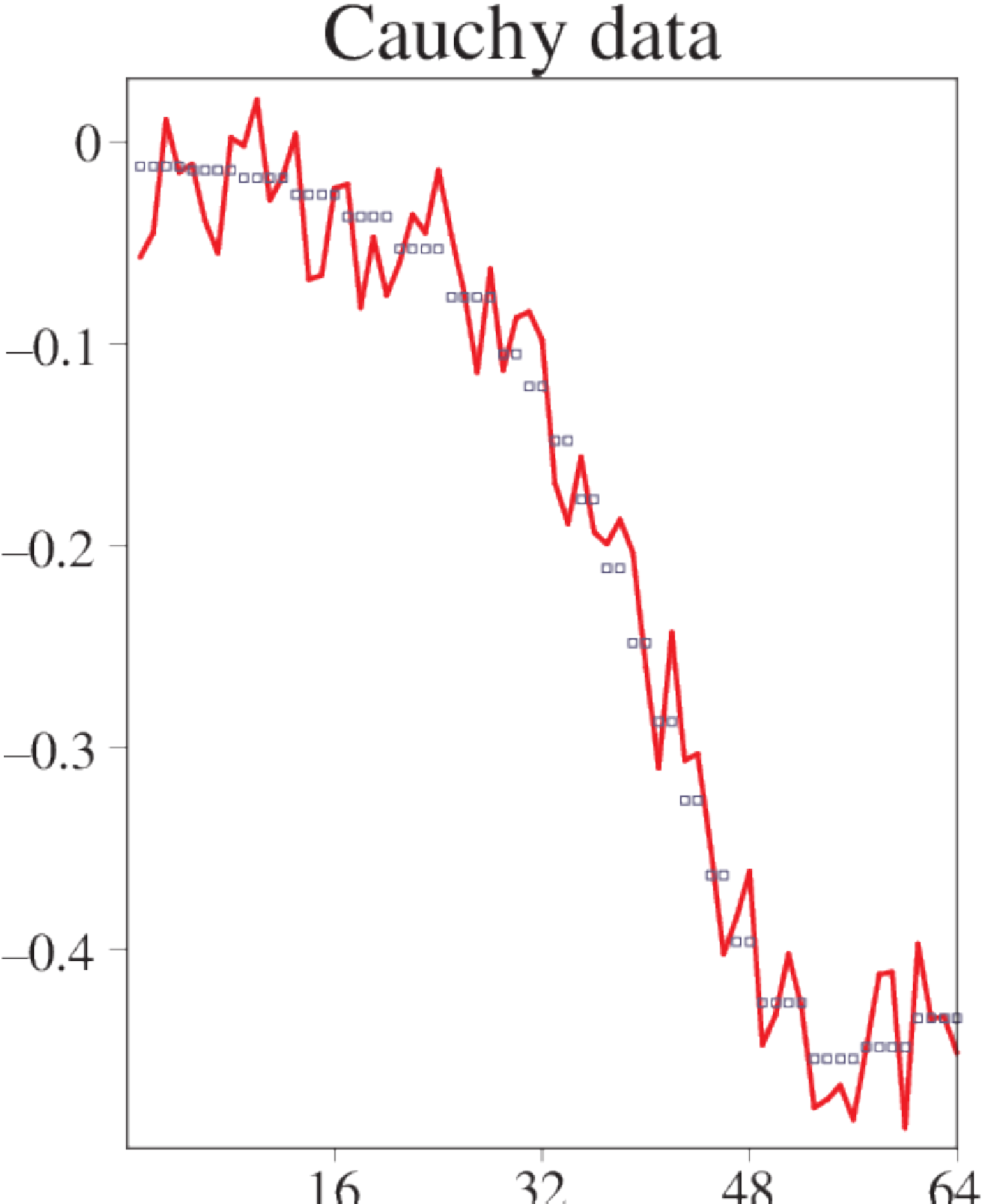} \hfil
             \epsfxsize4.2cm \epsfysize3.6cm
             \epsfbox{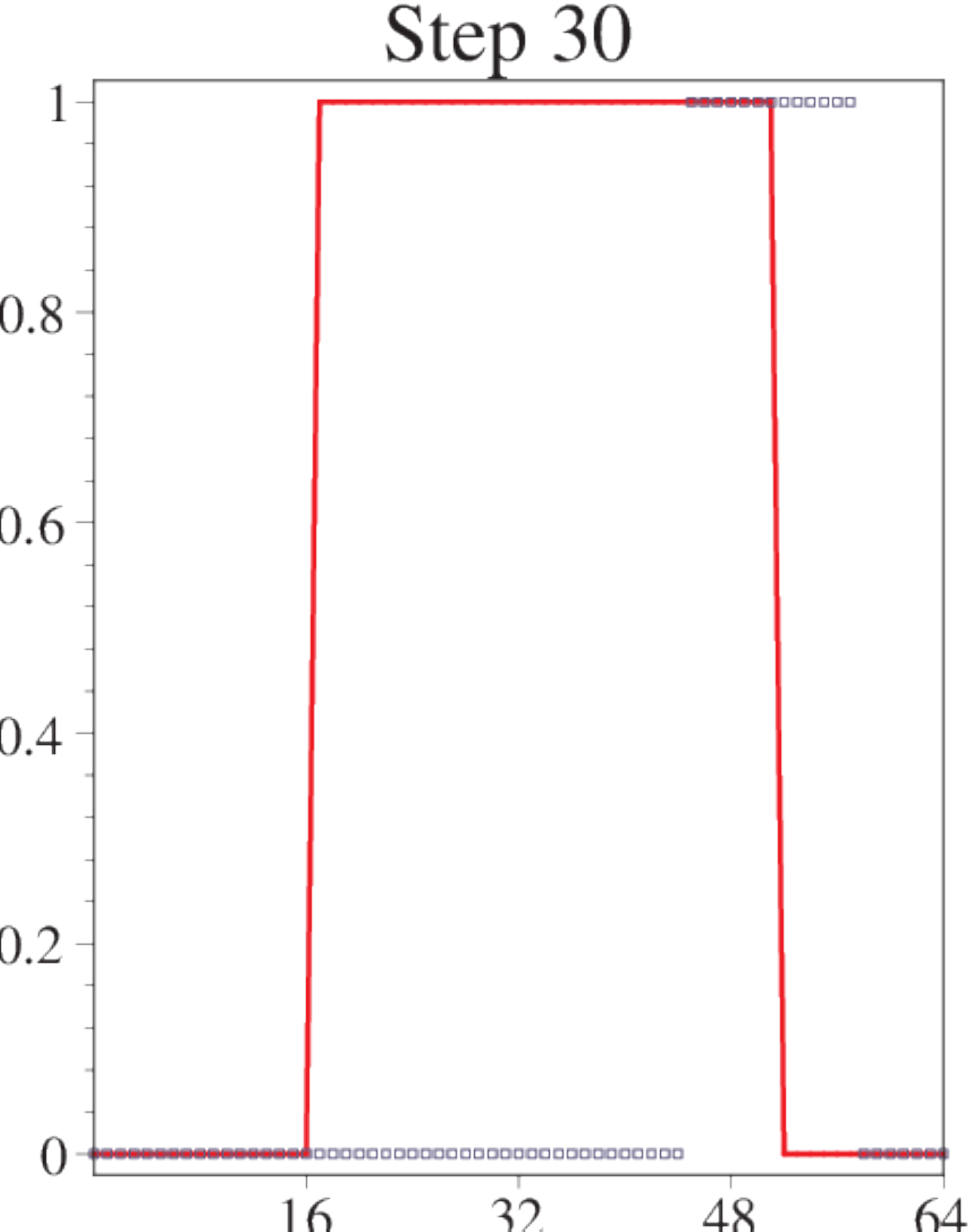} \hfil
             \epsfxsize4.2cm \epsfysize3.6cm
             \epsfbox{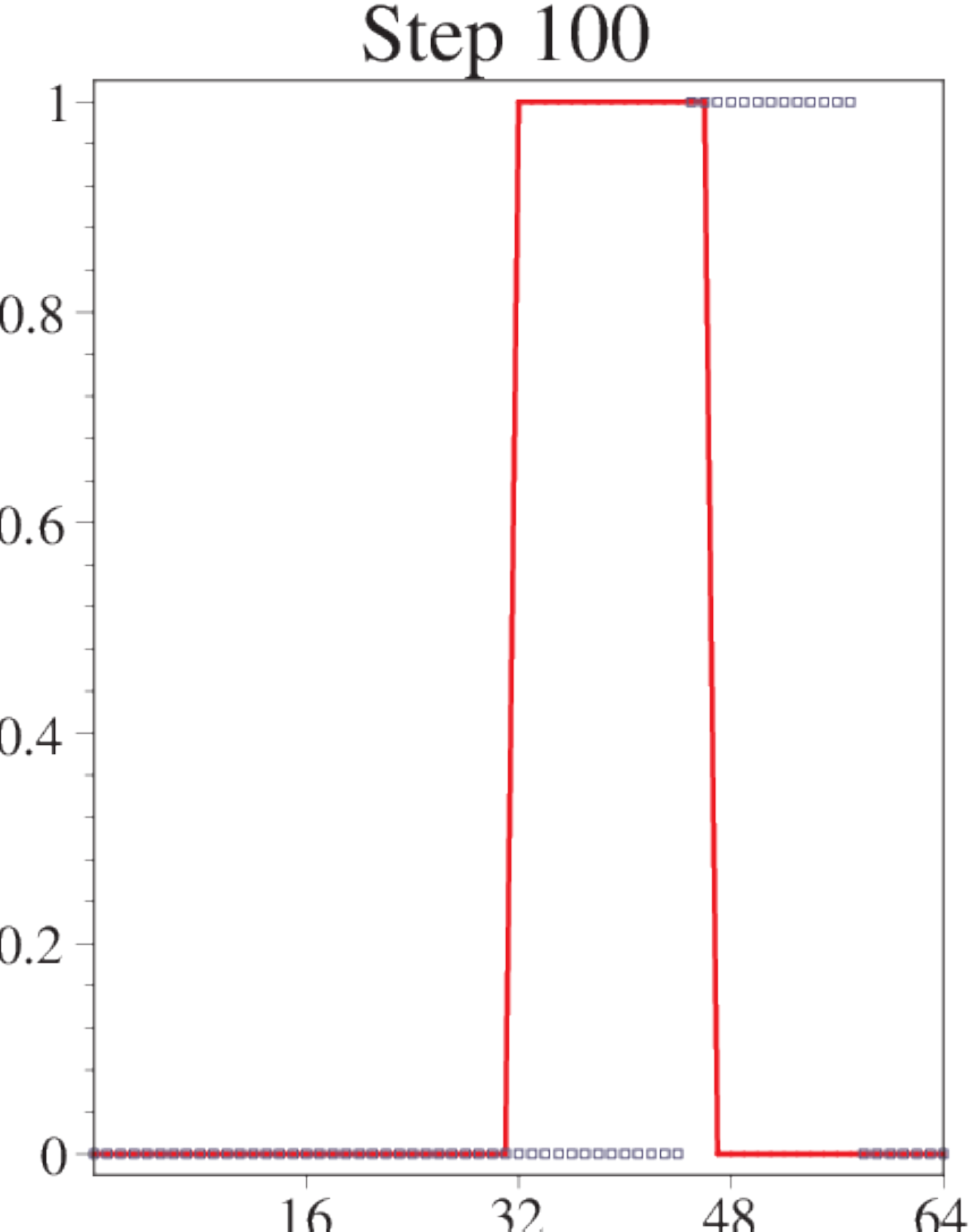} } \medskip
\centerline{ \epsfxsize4.2cm \epsfysize3.6cm
             \epsfbox{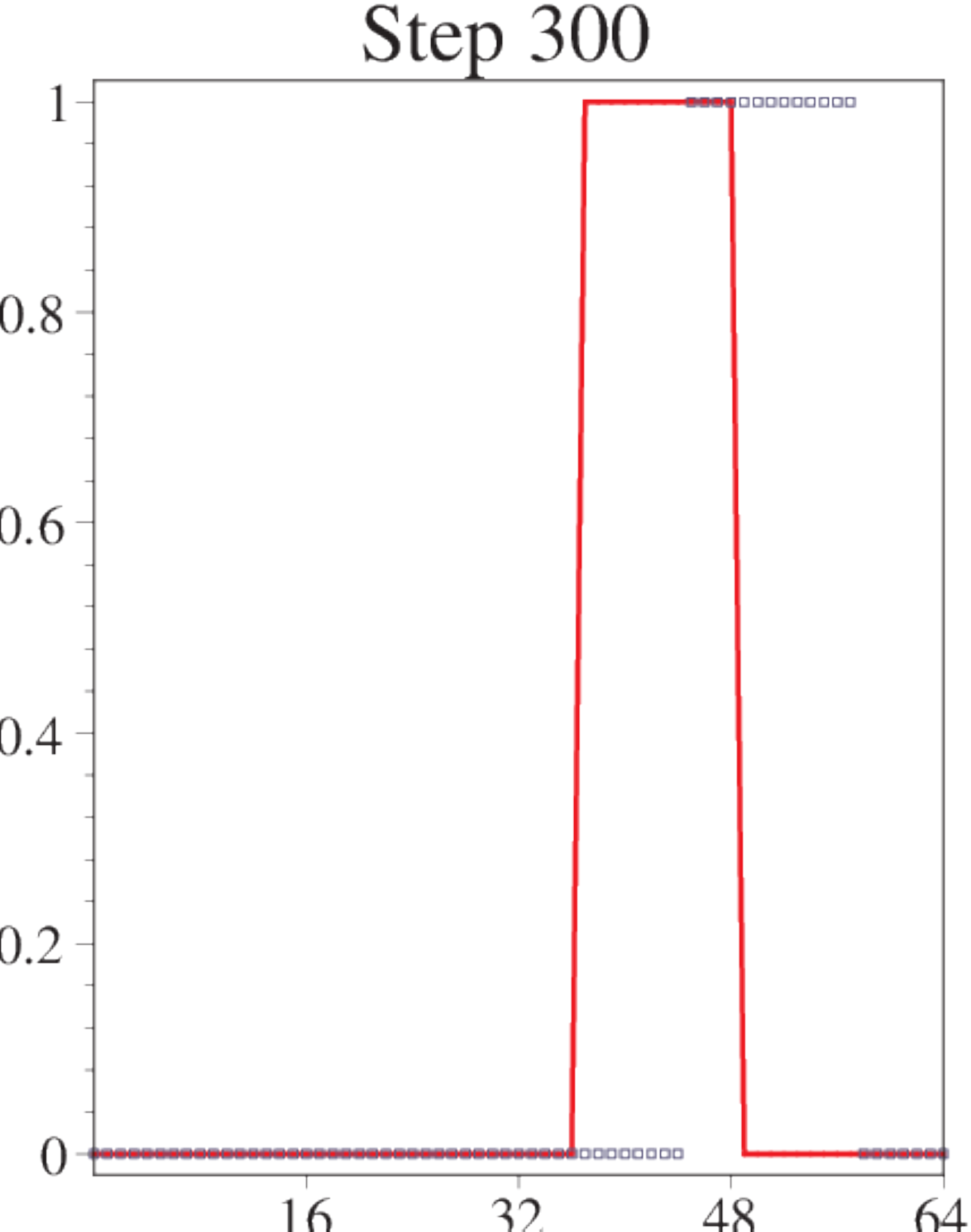} \hfil
             \epsfxsize4.2cm \epsfysize3.6cm
             \epsfbox{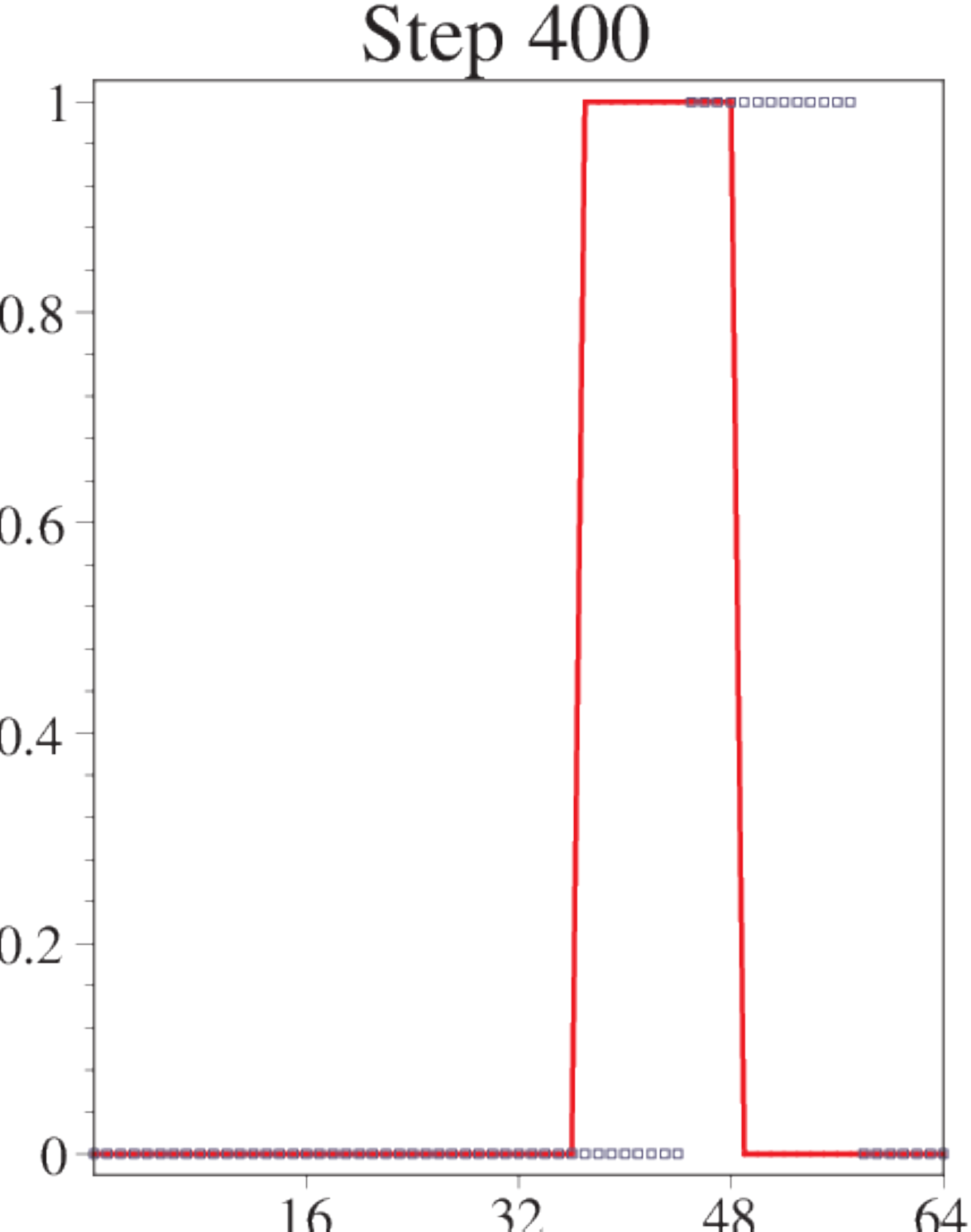} \hfil
             \epsfxsize4.2cm \epsfysize3.6cm
             \epsfbox{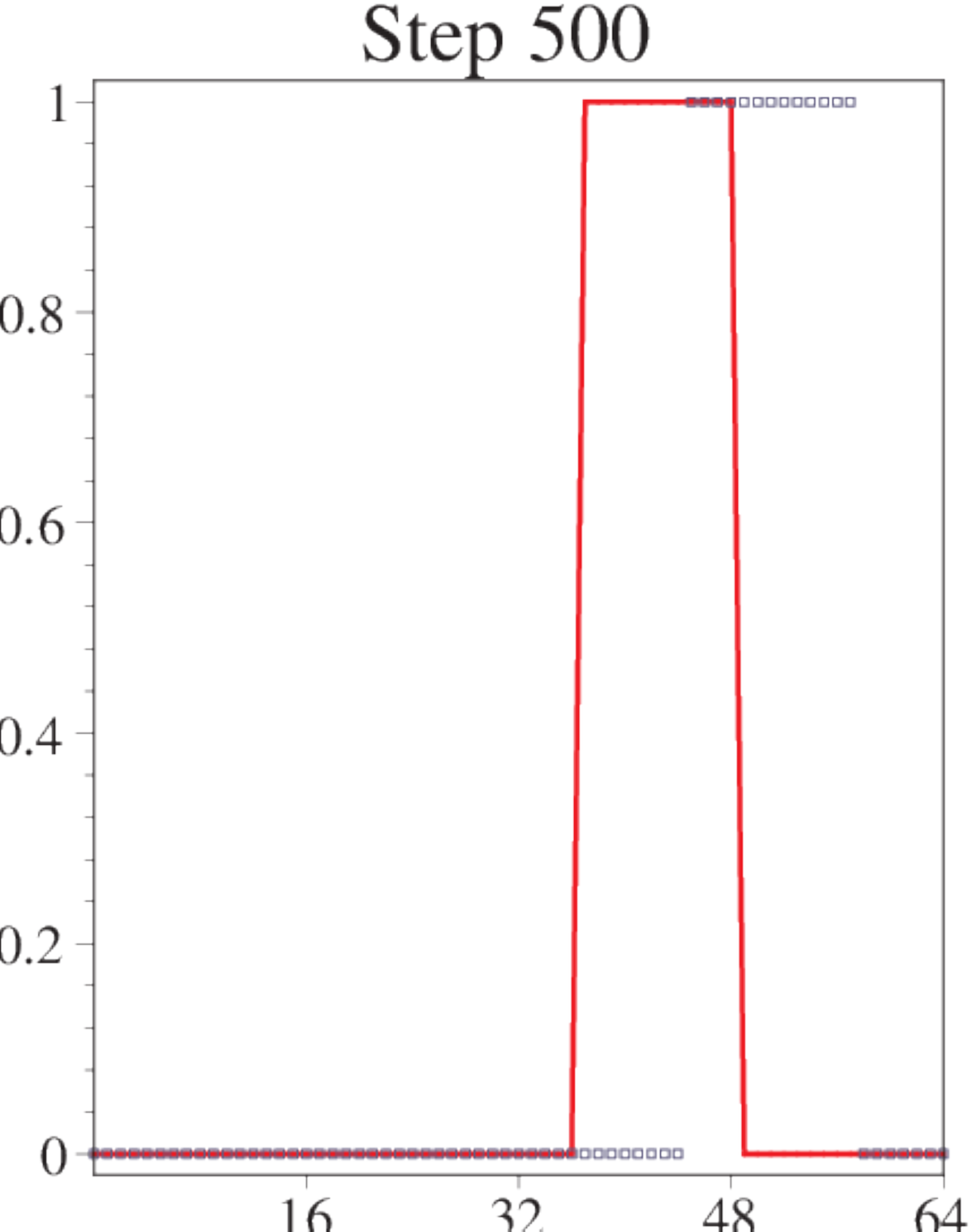} }
\caption{\footnotesize Third numerical experiment:
On the top left, the exact Cauchy data $g_{2,2}$ (dotted blue line) and
the data corrupted with 10\% random noise (solid red line).
The other pictures show the evolution of the level set method for problem
(CP$)_2$ with noisy data after 30, 200, 300, 400 and 500 iterative steps.
The solid (red) line represent the iteration and the dotted (blue) line the
exact solution. \label{fig:exp3-evol-ls}}
\end{figure}

For the next numerical experiment we consider once more problem (CP$)_2$
in Subsection~\ref{ssec:5-2}. This time however, we perturbed the Cauchy
data $g_{2,2}$ with 10\% random noise (Figure~\ref{fig:exp3-evol-ls}).
The initial guess for the level set method is the same used in the second
experiment (see Figure~\ref{fig:exp2-evol-ls2}).

The performance of the level set method for exact and corrupted data can
be compared in Figures~\ref{fig:exp3-evol-ls} and~\ref{fig:exp2-evol-ls2}.
Notice that for the noisy data the best possible approximation is obtained
after 300 steps. We iterated further (until 500 steps) and observed that,
although the level set function oscillates, the corresponding level sets
remain almost unchanged.

Since we assume the Cauchy data to satisfy $g_1 = 0$ in our experiments,
it seems natural not to introduce noise in this component of the data.
We conjecture this is the main reason for obtaining stable reconstruction
results for this exponentially ill-posed problem in the presence of high
levels of noise (in \cite{EL01} only 5\% noise is used; moreover it is
not white noise but corresponds to an eigenfunction of $L$).

\section{Conclusions} \label{sec:6}

In this article two possible level set approaches for solving elliptic
Cauchy problems are considered. For each one of them a corresponding
framework is established and a convergence analysis is provided
(monotony, convergence, stability results). Further we discuss the
numerical realization of the second level set approach.
Three different numerical experiments illustrate relevant features of this
level set method: rates of convergence, adaptability to identify
non-connected inclusions, robustness with respect to noise.

\section*{Acknowledgments}
The authors would like to thank Dr. Stefan Kindermann (Linz) for the
stimulating discussions during a research stay of the second author
at RICAM institute.
The work of M.M.A. is supported by the Brazilian National Research
Council CNPq.
A.L. acknowledges support from CNPq, grant 306020/2006-8.


\end{document}